\documentclass[12pt,reqno]{amsart}
\usepackage{amssymb,amsthm,amsfonts,amstext}
\usepackage{amsmath}
\usepackage{zito}
\usepackage{mathrsfs}
\usepackage{fourier}

\usepackage{dsfont}

\usepackage{xcolor}

\makeatletter \@addtoreset{equation}{section}  
\makeatother

\newtheorem{theorem}{Theorem}[section]
\newtheorem{lemma}[theorem]{Lemma}
\newtheorem{proposition}[theorem]{Proposition}

\theoremstyle{remark}
\newtheorem{remark}[theorem]{Remark}

\newcommand{\mc}[1]{{\mathcal #1}}

\newcommand{\bb}[1]{{\mathbb #1}}

\DeclareMathOperator{\Var}{Var}

\DeclareMathOperator{\Law}{Law}
\DeclareMathOperator{\Div}{div}
\DeclareMathOperator{\Cov}{Cov}
\DeclareMathOperator{\mix}{mix}
\DeclareMathOperator{\Beta}{Beta}

\newcommand{\sumxai}{{\sum_{x \in A^i}}}
\newcommand{\sumiq}{{\sum_{i \in Q}}}

\newcommand{\xis}{{\xi_{\hspace{-1.2pt}s}}}
\newcommand{\xiu}{{\xi_{\hspace{-1.2pt}u}}}

\newcommand{\zetat}{\zeta_{\hspace{-1.5pt}t}}
\newcommand{\zetas}{\zeta_{\hspace{-1.5pt}s}}

\newcommand{\Gamman}{{\Gamma_{\!n}}}
\newcommand{\munab}{{\mu_{\!a,b}^n}}
\newcommand{\xit}{\xi_{\hspace{-1.2pt}t}}

\newcommand{\vecX}{{\vec{X}}}

\newcommand{\va}{{\vec{a}}}
\newcommand{\vecM}{{\vec{M}}}
\newcommand{\vecMt}{{\vec{M}_t}}
\newcommand{\vecF}{{\vec{F}}}
\newcommand{\vecG}{{\vec{G}}}

\newcommand{\vecm}{{\vec{m}}}
\newcommand{\vecmt}{{\vec{m}_t}}
\newcommand{\vecms}{{\vec{m}_s}}

\newcommand{\vecmzero}{{\vec{m}_0}}

\newcommand{\vecxi}{{\vec{\xi}}}

\newcommand{\vecxit}{{\vec{\xit}}}
\newcommand{\vecxis}{{\vec{\xis}}}
\newcommand{\vecxiu}{{\vec{\xiu}}}

\newcommand{\veczetat}{{\vec{\zetat}}}
\newcommand{\veczetas}{{\vec{\zetas}}}
\newcommand{\veczetazero}{\vec{\zeta}_{\!0}}

\newcommand{\vecyt}{{\vec{y}_t}}
\newcommand{\vecys}{{\vec{y}_s}}
\newcommand{\vecyzero}{{\vec{y}_0}}
\newcommand{\vecYt}{{\vec{Y}_t}}
\newcommand{\vecYtn}{{\vec{Y}_{t_n}}}

\newcommand{\ytast}{{y_{\!t}^\ast}}
\newcommand{\ztast}{{z_{\!t}^\ast}}

\newcommand{\veczt}{{\vec{z}_t}}

\newcommand{\veczzero}{{\vec{z}_0}}
\newcommand{\vecZt}{{\vec{Z}_t}}
\newcommand{\vecZtn}{{\vec{Z}_{t_n}}}

\newcommand{\vecx}{{\vec x}}

\newcommand{\Ln}{{\Lambda_n}}
\newcommand{\On}{{\Omega_n}}
\newcommand{\Onast}{{\Omega_n^\ast}}
\newcommand{\sumxln}{{\sum_{x \in \Ln}}}

\newcommand{\pens}{{\bb P_{\!n}^\sigma}}
\newcommand{\bens}{{\bb E_n^\sigma}}
\newcommand{\penus}{{\bb P_{\!n}^{u(\sigma)}}}
\newcommand{\vq}{{{\bb V}_Q}}

\title[Mixing of noisy voter model]{Thermalization and convergence to equilibrium of the noisy voter model}

\author{Enzo Aljovin}
\address[E. Aljovin]{Instituto de Matem\'atica Pura e Aplicada, Estrada Dona Castorina 110, 22460-320
	Rio de Janeiro, Brasil}
\email{enzo.aljovin@gmail.com}

\author{Milton Jara}
\address[M. Jara]{Instituto de Matem\'atica Pura e Aplicada, Estrada Dona Castorina 110, 22460-320
	Rio de Janeiro, Brasil}
\email{mjara@impa.br}

\author{Yangrui Xiang}
\address[Y. Xiang]{Instituto de Matem\'atica Pura e Aplicada, Estrada Dona Castorina 110, 22460-320
	Rio de Janeiro, Brasil}
\email{yangrui94.xiang@gmail.com}

\begin{document}

\begin{abstract}
We investigate the convergence towards equilibrium of the noisy voter model, evolving in the complete graph with $n$ vertices. The noisy voter model is a version of the voter model, on which individuals change their opinions randomly due to external noise. Specifically, we determine the profile of convergence, in Kantorovich distance (also known as $1$-Wasserstein distance), which corresponds to the Kantorovich distance between the marginals of a Wright-Fisher diffusion and its stationary measure. In particular, we demonstrate that the model does not exhibit cut-off under natural noise intensity conditions. In addition, we study the time the model needs to forget the initial location of particles, which we interpret as the Kantorovich distance between the laws of the model with particles in fixed initial positions and in positions chosen uniformly at random. We call this process \emph{thermalization} and we show that thermalization does exhibit a cut-off profile. Our approach relies on Stein's method and analytical tools from PDE theory, which may be of independent interest for the quantitative study of observables of Markov chains.
\end{abstract}

\maketitle

\section{Introduction}

\subsection{Overview} The noisy voter model is an irreducible, continuous-time Markov chain and it can be regarded as a perturbation of the usual Voter Model. The classical theory of Markov chains guarantees the existence of a unique invariant (probability) measure. Nevertheless, the exponential convergence rate to the stationary state, which is determined by the \emph{spectral gap} of the chain, is not always sufficient for practical purposes. Recent developments in Markov chain Monte Carlo (MCMC) algorithms have highlighted the need for deeper understanding on the times the chain takes to get within a given threshold of the stationary state, the so-called \emph{mixing time}. Furthermore, finding computationally tractable modifications to existing Markov chains, that speed up their mixing time, is an ongoing challenge. The potential benefits of such challenges include simulation of near-equilibrium samples, in a computationally feasible way.

The goal of this paper is to investigate the mixing properties of the noisy voter model on the complete graph, on certain noise intensity regimes, taking the population size as a scaling parameter. Moreover, we explore the influence of the initial state of the chain on its mixing properties.

The main contribution of this article is to quantitatively determine the \emph{profile of convergence} towards the equilibrium measure. In other words, we prove that the Kantorovich distance between the law of the chain at time $nt$ and its stationary measure, converges to a \emph{continuous profile}. This profile arises as the profile of convergence of a Wright-Fisher diffusion, which can be regarded as the limit of a discrete model for genetic drift in a population where immigration takes place. Furthermore, the continuity of the profile allows us to conclude that the convergence occurs \emph{smoothly}, in the sense that the family of chains does not display cut-off.

As the proof of the result above shows, the convergence to equilibrium is driven by the convergence of the \emph{empirical density of particles} to the solution of the Wright-Fisher diffusion. It is natural therefore to ask how other observables of the model equilibrate. We answer this question by studying the dependence of the law of the model, \emph{given the initial number of particles}, on the initial position of those particles. We show that the process forgets the initial position of particles in times of order $\frac{1}{2} \log n$, much shorter than the times of order $nt$, at which the density of particles converge to the Wright-Fisher diffusion. If we interpret the number of particles as \emph{wells} in the sense of metastability, this corresponds to what is known in the literature as \emph{thermalization}. This thermalization does satisfy cut-off with an explicit cut-off profile.

\subsection{Noisy voter model}
In the decade of the 70s, Clifford and Sudbury \cite{CliSud} and Holley and Liggett \cite{HolLig} introduced the \emph{Voter Model}, an interacting particle system that serves as a model for a variety of social, biological and physical phenomena. The Voter Model can describe opinion dynamics driven by social interactions in a network, competition for territory between two species \cite{CliSud}, neutral genetic drift in a population \cite{KimOht}, among others. In the Voter Model, each individual has its own opinion and they interact copying opinion of others, at random times. It is worth mentioning that the Voter Model has absorbing states, which are consensus. In the case of finite population, consensus is reached in finite time and no further evolution is observed afterwards. 

Later, Granowsky and Madras \cite{GraMad} introduced an ergodic variant of the Voter Model, which is nowadays known as the \emph{noisy voter model}. The noisy voter model is a superposition of the dynamics of the Voter Model and a re-randomization of the opinions. Specifically, at rate $a$ (resp.~$b$), one vertex is selected uniformly at random and its opinion is set to $1$ (resp.~$0$). In the context of genetics, the rates of re-randomization can be rephrased as \emph{mutation rates}. As the name suggests, re-randomization drives the system away from consensus. Thus, two major effects compete in the noisy voter model. On one hand, the voter dynamics drives the system towards consensus. While, on the other hand, spontaneous changes of opinion promote coexistence of opinions. Additionally, since the chain is ergodic, it has a unique non-trivial invariant measure. In the case of the noisy voter model evolving on the complete graph, the invariant measure corresponds can be obtained by choosing the number of particles with Beta-Binomial law with parameters $(n,a,b)$ and then choosing the positions of such particles uniformly at random, where $n$ corresponds to the size of the underlying graph and $a$, $b$ are the noise intensity parameters.

\subsection{Mixing properties} 
It is well established that if the size of the state space of an irreducible Markov chain is fixed, convergence towards the equilibrium measure occurs exponentially fast \cite[Chapter 4]{LevPer}. This result is sometimes called the \emph{Fundamental Theorem of Markov Chains}. However, in this paper we adopt a different point of view, pioneered by Aldous and Diaconis \cite{AldDia}. The core of this approach is to consider a parametrized family of Markov chains (being the size of the state space the most common parameter, see however \cite{BarJar}), and then to fix a threshold distance to stationarity. The arising question is how much time does the chain needs to reach that proximity to its stationary state, with time scaled as a function of the size of the state space. One of the main differences with the classical convergence theorem approach, is that the associated constants are size-dependent and they usually grow at a faster rate, which makes this approach unsuitable for practical applications. 

In order to be more precise, let us assume that the parameter $n \in \bb N$ and let $(\Omega_n, d_n)$ be the state space of the chain of parameter $n$, equipped with a distance $d_n$ in the space of measure on $\Omega_n$. Let $(P_t^n; t \geq 0)$ be the semigroup of the chain and let $\pi^n$ be its invariant measure. Given $\varepsilon >0$, the $\varepsilon$-\emph{mixing time} of the chain in $(\Omega_n, d_n)$ is defined as
\[
t_{\mix}^n(\varepsilon) := \inf\big\{ t \geq 0; \sup_{\eta \in \On} d_n \big(P_t^n(\eta,\cdot), \pi^n \big) \leq \varepsilon \big\}.
\]
The most common choice for $d_n$ is the \emph{total variation distance}, which has the advantage of being well defined regardless of the state space $\On$. From a probabilistic point of view, it is the most natural choice, due to a very useful duality relation. Nonetheless, the definition of mixing time takes as initial condition the \emph{worst possible case}. This choice makes sense from a purely theoretical point of view, where definitions should not take into account any particularity of the underlying space. From the point of view of real-world applications, it makes perfect sense to incorporate the intuition we have about the model we are interested about. For example, it is of interest to initialize an MCMC algorithm from the most convenient initial condition. In this paper, we aim to investigate the dependence of the mixing properties with respect to the initial state of the chain. 

While mixing times provide valuable information regarding \emph{when} convergence towards equilibrium occurs, they do not fully capture \emph{how} the convergence actually happens. In the early works of Aldous, Diaconis and coauthors, methods from representation theory allowed them to access the \emph{profile of convergence} in the cut-off phenomenon. More precisely, they showed that there exist sequences $t^n_{\mix}$, $\omega^n_{\mix}$ such that $\omega_{\mix}^n \ll t_{\mix}^n$ and 
\[
\lim_{n \to \infty} d_n \big(P_{t_{\mix}^n + q \omega_{\mix}^n}^n(\eta,\cdot), \pi^n \big) = \mc G(q),
\]
where $\mc G$ is a non-trivial function called the \emph{convergence profile}. See \cite{AldStras} for a review and further references, and see in particular Example (3.19) of \cite{AldStras}. In a series of papers \cite{BarJar, BarJar2, BarPar, BarHoePar}, the authors describe the cut-off profile for superstable dynamical systems perturbed by small Brownian and/or L\'evy noise.

In a series of breakthrough papers \cite{Lac1,Lac2,Lac3}, the author studied the mixing properties of the exclusion process on the interval and on the circle. The author obtained sharp estimates on the mixing time and proved that these models exhibit the cut-off phenomenon, with a \emph{universal} profile function. In \cite{GonJarMarMen}, the dependence of this cut-off profile with respect to non-equilibrium measures that commonly appear in the theory of hydrodynamic limits is investigated.

Choosing $a = pn$, $b=(1-p) n$ with $p \in (0,1)$, our noisy voter model is one of the models considered in \cite{CoxPerSte}. In this paper, the authors show that the noisy voter model has cut-off at time $\frac{1}{2} \log n$, with cut-off window of order $\mc O(1)$. Observe that in this situation, the number of particles  has a Beta-Binomial law of parameters $(n, pn, (1-p)n)$ and in particular it has Gaussian fluctuations of size $\mc O(\sqrt n)$ around its mean $pn$. In our situation, the number of particles has a law extended up to the whole interval $\{0,1, \dots, n\}$. This fact changes completely the behavior of the model. We will see that in this case the convergence is not abrupt, and it happens smoothly in a time scale $\mc O(n)$. Nevertheless, we will see that at times of order $\mc O(\frac{1}{2} \log n)$, the model forgets the initial \emph{positions} of the particles, and it only remembers the initial \emph{number} of particles. This matches the predictions of \cite{CoxPerSte}. In order to bridge our results with the results of \cite{CoxPerSte}, we formulate the following conjecture. Let $(r_n; n \in \bb N)$ be such that
\[
\lim_{n \to \infty} r_n = +\infty \text{ and } \lim_{n \to \infty} \frac{r_n}{n} = 0.
\]
Then the noisy voter model with parameters $a r_n$, $b r_n$ has Gaussian profile cut-off at time $\frac{r_n}{2n(a+b)} \log r_n$ with window $\frac{r_n}{n(a+b)}$.

\subsection{Main ideas}
Let us describe the main ideas of the proofs presented in this paper. At scales of order $t \ll n^{2/9}$, we show a quantitative Central Limit Theorem (CLT) for the fluctuations of the density of particles. In order to this, we approximate the vector of \emph{local densities} by the solution of a time-inhomogeneous Ornstein-Uhlenbeck equation, in the Kantorovich distance. Such a process is Gaussian, from where an asymptotic analysis allows to compare the system starting from different initial conditions. The quantitative CLT is proved by means of \emph{Stein's method}, in particular the semigroup method introduced by Barbour \cite{Bar}. At the larger time scale $t= \mc O(n)$, we show that the one-time marginals of the global empirical density of particles can be quantitatively approximated in Kantorovich distance by the solution of the Wright-Fisher equation. Again, the semigroup approach to Stein's method is the main technical tool. We thanks the authors of \cite{HerJar} for providing an early manuscript, where this strategy was introduced in the context of the stochastic Curie-Weiss process, and for allow us to use Lemma \ref{reglem}.

\section{Notation and results}

\label{s2} 
\subsection{The noisy voter model} Let $n \in \bb N$ be a scaling parameter and let $\Ln:= \{1,\dots,n\}$, $Q :=\{0,1\}$ and $\On:= Q^{\Ln}$. For $\eta := ( \eta_x;  x \in \Ln) \in  \On$ and $x \in \Ln$, let $\eta^x \in \On$ be given by
\[
\eta_z^x:=
\left\{
\begin{array}{c@{\;;\;}l}
1-\eta_x & z =x,\\
\eta_z & z \neq x.
\end{array}
\right.
\]
The elements $x \in \Ln$ are called \emph{sites} and the elements $\eta \in \On$ are called \emph{particle configurations}. We say that a configuration $\eta \in \On$ has a \emph{particle} at site $x \in \Ln$ if $\eta_x =1$. Otherwise, we say that the site $x$ is \emph{empty}. 

Let $X : \On \to [0,\infty)$ be given by
\[
X(\eta) := \!\!\sumxln \!\eta_x
\]
for every $\eta \in \On$, that is, $X(\eta)$ is the \emph{number of particles} of the configuration $\eta \in \On$. Let $a,b >0$ be fixed and for every $x \in \Ln$, let $c_x: \On \to [0,\infty)$ be given by
\[
c_x(\eta) := (1-\eta_x) \big(a+X(\eta)\big) + \eta_x \big(b+n-X(\eta)\big)
\]
for every $\eta \in \On$. For $f: \On \to \bb R$, let $L_n f: \On \to \bb R$ be given by
\[
L_n f(\eta) := \frac{1}{n} \sumxln c_x(\eta) \big( f(\eta^x) - f(\eta) \big)
\]
for every $\eta \in \On$. The linear operator $L_n$ defined in this way turns out to be the generator of a continuous-time Markov chain $(\eta(t); t \geq 0 )$ that we call the \emph{noisy voter model} in the complete graph $\Ln$. 

It will be useful to define the \emph{carr\'e du champ} $\Gamman$ associated to the operator $L_n$. The bilinear operator $\Gamman$ is defined as 
\[
\Gamman (f,g) := L_n(fg) - f L_n g - g L_n f
\]
for every $f,g: \On \to \bb R$. As usual, we will write $\Gamman f$ instead of $\Gamman (f,f)$. Observe that
\[
\Gamman f(\eta) = \frac{1}{n} \sumxln c_x(\eta) \big( f(\eta^x) - f(\eta) \big)^2 
\]
for every $\eta \in \On$. 

The noisy voter model $(\eta(t); t \geq 0)$ is an irreducible Markov chain on a finite state space, and therefore it has a unique invariant measure $\munab$ with total support in $\On$. Moreover, the measure $\munab$ is reversible and explicit. For $\alpha, \beta >0$, let
\[
B(\alpha, \beta) := \int_0^1 t^{\alpha-1} (1-t)^{\beta-1} dt
\]
be the \emph{Beta function}. It can be verified that
\[
\munab(\eta) = \frac{B(a+X(\eta), b+n-X(\eta))}{B(a,b)}
\]
for every $\eta \in \On$. Observe that $\munab$ depends on $\eta$ only through the number of particles $X$. Therefore, the measure $\munab$ is uniform on each of the sets
\[
\Omega_{n,\ell} := \big\{\eta \in \On; X(\eta) = \ell\big\}, \ell =0,1,\dots,n.
\]

We will denote by $\bb P_{\!n}^\mu$ the law of the chain $(\eta(t); t \geq 0)$ with initial measure $\mu$ and we will denote by $\bb E_n^\mu$ the expectation with respect to $\bb P_{\!n}^\mu$. If the initial measure $\mu$ is a delta of Dirac supported at $\sigma \in \On$, we will write $\pens$ and $\bens$ instead of $\bb P_{\!n}^{\delta_\sigma}$ and $\bb E_n^{\delta_\sigma}$.

The main objective of this article is the study of the convergence of the law of $\eta(t)$ to $\munab$ as a function of the initial condition $\sigma \in \On$ of the chain $(\eta(t); t \geq 0)$. Our main tool will be the derivation of \emph{quantitative} central limit theorems (QCLT) for various quantities of interest. In order to state these QCLTs, we need to introduce a distance on the space of probability measures.

\subsection{The Kantorovich distance} Let $(\mc E,d)$ be complete, separable metric space. Let $\mc B(\mc E)$ be the Borel $\sigma$-algebra associated to $(\mc E,d)$. For $f: \mc E \to \bb R$, define
\[
\llbracket f \rrbracket := \sup_{x \neq y} \frac{|f(y) - f(x) |}{d(x,y)}.
\] 
If $\llbracket f \rrbracket <+\infty$, we say that $f$ is \emph{Lipschitz}. Let $x_0 \in \mc E$ be fixed. We say that a probability measure $\mu$ in $(\mc E, \mc B(\mc E))$ \emph{belongs to} $\mc P_1(\mc E)$ if
\[
\int d(x_0,x) \mu(dx) <+\infty.
\]
Observe that in that case, for every $y \in \mc E$,
\[
\int d(y,x) \mu(dx) \leq d(y,x_0) + \int d(x_0,x) \mu(dx) <+\infty
\]
and in particular the definition of $\mc P_1(\mc E)$ does not depend on the choice of $x_0$. For $\mu, \nu \in  \mc P_1(\mc E)$, we define the \emph{Kantorovich distance}\footnote{This distance is sometimes called in the literature the \emph{$1$-Wasserstein distance.}} between $\mu$ and $\nu$ as
\[
d_K(\mu,\nu) := \sup_{\llbracket f \rrbracket \leq 1} \Big| \int f d\mu - \int f d\nu\Big|.
\]
The space $\mc P_1(\mc E)$ turns out to be complete and separable with respect to the Kantorovich distance. The relevance of the Kantorovich distance in probability theory is given by the so-called \emph{Kantorovich-Rubinstein duality}, which now we state. We say that a measure $\pi$ in the product space $\mc E \times \mc E$ is a \emph{coupling} of $\mu$ and $\nu$ if
\[
\pi(A \times \mc E) = \mu(A) \text{ and } \pi( \mc E \times A) = \nu(A)
\]
for every open set $A \subseteq \mc E$. Let $\Pi_{\mu,\nu}$ be the set of couplings of $\mu$ and $\nu$. The Kantorovich-Rubinstein duality states that
\begin{equation}
\label{dual}
d_K(\mu, \nu) = \inf_{\pi \in \Pi_{\mu,\nu}} \iint d(x,y) \pi(dx dy)
\end{equation}
and that the infimum is attained at least in one coupling $\pi_\ast$. 

For $\mc E$-valued random variables $X$, $Y$, we use the notation 
\[
d_K(X,Y) := d_K(\mu_X, \mu_Y),
\]
where $\mu_X$ and $\mu_Y$ are the laws of the variables $X$, $Y$ in $(\mc E, \mc B(\mc E))$. If there is no risk of confusion, we will also write $d_K(X,\nu)$ instead of $d_K(\mu_X,\nu)$.

In this paper, we will consider the Kantorovich distance in the spaces $\bb R$, $\bb R^Q$ and $[0,1]$ equipped with the Euclidean distance, and in the spaces $\On$ equipped with the \emph{Hamming distance} defined as
\[
\bb H_n(\eta, \sigma) := b_n \sumxln |\eta_x - \sigma_x|.
\]
The \emph{scale} $b_n$ is case dependent; depending on the situation, we will use the scales $b_n = \frac{1}{n}$ and $b_n = \frac{1}{\sqrt n}$.

Sometimes we will need to compare the law of the same random variable with respect to different measures. In order to do this, given a space $\Omega$, a measure $\mu$ in $\Omega$ and a random variable $X: \Omega \to \mc E$, we will denote  by $\Law(X; \mu)$ the law of $X$ in $(\mc E, \mc B(\mc E))$ with respect to $\mu$.

\subsection{Mixing times and sharp convergence}

In the study of the evolution of the noisy voter model with parameters $a,b$, we split the dynamics in two conceptually different phases. First we loot at times of order up to $\mc O(\frac{1}{2} \log n)$. During this stage of the dynamics,  we will show that the dynamics remembers the initial number of particles, and it starts to forget the initial positions of those particles. More precisely, let $\bar{\sigma}^n, \tilde{\sigma}^n \in \Omega_n$ be such that
\[
X(\bar{\sigma}^n) = X(\tilde{\sigma}^n) 
\]
for every $n \in \bb N$, and 
\[
\lim_{n \to \infty} \frac{X(\bar{\sigma}^n)}{n} = m_0 \in (0,1).
\]
We will prove that for every $\tau \in \bb R$, setting
\[
t_n = := \frac{1}{2} \log n+ \log m_0(1-m_0),
\]
\[
\lim_{n \to \infty} \frac{1}{\sqrt n} d_K \big( \Law \big(\sigma(t_n + \tau); \bb P_{\!n}^{\bar{\sigma}^n} \big), \Law \big(\sigma(t_n + \tau); \bb P_{\!n}^{\tilde{\sigma}^n}\big) \big) = e^{-\tau}.
\]
This is a consequence of Theorem \ref{phase1} below, which states a more quantitative version of this result. The scaling factor $\frac{1}{\sqrt n}$ in the Hamming distance turns to be the natural one, since the difference between the local density of particles and the total density of particles has fluctuations of order $\mc O(\frac{1}{\sqrt n})$. 

At times of order $\mc O(n)$, the model starts to forget the initial number of particles. As mentioned above, at this time scale the number of particles can be approximated by the solution of the Wright-Fisher equation. The natural scaling factor in the Hamming distance is now $\frac{1}{n}$, and under the assumptions above, we can show that
\begin{equation}
\label{tacna}
\lim_{n \to \infty} \frac{1}{n} d_K\big(\Law \big(\sigma(nt); \bb P_{\!n}^{\bar{\sigma}^n} \big), x_t(m_0) \big) =0,
\end{equation}
where $(x_t(m_0), t \geq 0)$ is the solution of the Wright-Fisher equation of parameters $a, b$ and initial condition $m_0$. This is a consequence of Theorems \ref{t2} and \ref{phase1} below.

In order to analyze the dependence of the mixing time on the initial condition, let us define, for $n \in \bb N$, $\varepsilon > 0$ and $\sigma \in \On$,
\[
t_{\mix}^n(\varepsilon; \sigma) := \inf \big\{ t \geq 0; \frac{1}{n} d_K \big(P_t^n(\sigma, \cdot ), \munab \big) \leq \varepsilon \big\}.
\]
Let $x_\infty$ be a random variable with law  $\Beta(a,b)$. Using \eqref{tacna} and the fact that the function
\[
\mc G(t; m):= d_K(x_t(m); x_\infty)
\]
is continuous and strictly decreasing, we see that
\begin{equation}
\label{ilo}
\lim_{n \to \infty} \frac{t_{\mix}^n(\varepsilon; \bar{\sigma}^n)}{n} = \mc G^{-1}(\varepsilon; m_0),
\end{equation}
which shows that the noisy voter model does not have cut-off in Kantorovich distance.

\begin{remark}
It is possible to prove that Theorem \ref{phase1} also holds in total variation distance. Relation \eqref{ilo} also holds if we define $\mc G$ and $t_{\mix}^n$ using total variation distance instead of Kantorovich distance. Since the proofs are extremely technical, they use the results in Kantorovich distance as input and the paper is already long enough, we decided to state our results only in terms of Kantorovich distance.
\end{remark}

\section{The convergence of the density of particles}

\subsection{The QCLT for the number of particles}
In this section we will prove a QCLT for the density of particles. More precisely, let $M: \On \to [0,1]$ the density of particles defined as
\[
M(\eta) := \frac{1}{n} \sumxln \eta_x 
\] 
for every $\eta \in \On$. Observe that $M = \frac{X}{n}$. Let $(M_t; t \geq 0)$ be given by $M_t := M(\eta(t))$ for every $t \geq 0$. For every $m \in [0,1]$, let $(x_t(m); t \geq 0)$ be the solution of the \emph{Wright-Fisher diffusion equation} 
\begin{equation}
\label{wrifish}
d x_t = \big( a(1-x_t) - b x_t \big) dt + \sqrt{ 2 x_t(1-x_t)} d B_t
\end{equation}
with initial condition $x_0 = m$. Here $(B_t; t \geq 0)$  is a standard Brownian motion in $\bb R$. We have the following result:

\begin{theorem}
\label{t2}
There exists $C= C(a,b)$ finite such that
\[
d_K \big(M_{nt}, x_t(M_0) \big) \leq \frac{C \min\{1,t^{1/4}\}}{\sqrt n}
\]
for every $n \in \bb N$, every $t \geq 0$ and every initial condition $\sigma \in \On$.
\end{theorem}

\begin{proof}
The proof of this theorem uses three lemmas that will be stated and proved in subsequent sections. The idea is to use a parabolic version of \emph{Stein's method}. By the definition of the Kantorovich distance, our aim is to estimate the difference
\begin{equation}
\label{arica}
\big| \bens \big[ f(M_{nt}) \big] - \bb E \big[ f(x_t(M_0)) \big] \big|
\end{equation}
in terms of $\|Df \|_\infty$. Alternatively, one can estimate this difference in terms of the norms of higher-order derivatives of $f$ and rely on Lemma \ref{reglem}.
Let $\Lambda_{a,b}$ be the generator of the Wright-Fisher diffusion, that is, 
\[
\Lambda _{a,b} f(x) := \big( a(1-x) -b x \big) f'(x) + x(1-x) f''(x)
\]
for every $f \in \mc C_b^2([0,1]; \bb R)$ and every $x \in [0,1]$. For $f \in \mc C_b([0,1]; \bb R)$ abd $t >0$, let $(g_{s,t}; s \in [0,t])$ be the solution of the \emph{backwards Fokker-Planck equation}
\begin{equation}
\label{fokpla}
\big( \partial_{\!s} + \Lambda_{a,b} \big) g_s =0
\end{equation}
with final condition $g_{t,t} =f$. By It\^o's formula,
\[
\bb E[ f(x_t(m))] = g_{0,t} (m)
\]
for every $m \in [0,1]$. By Dynkin's formula applied to the function
\[
(s,\eta) \mapsto g_{s,t} \big(M(\eta) \big),
\]
we see that
\[
\bens \big[ f(M_{tn}) \big] = g_{0,t}( M_0) + \int_0^t \bens \big[ (\partial_{\!s} + nL_n) g_{s,t}(M_{ns}) \big] ds.
\]
Therefore, in order to estimate \eqref{arica}, it is enough to estimate the integral
\begin{equation}
\label{caldera}
\Big| \int_0^t  \bens \big[ (\partial_{\!s} + nL_n) g_{s,t}(M_{ns}) \big] ds \Big|.
\end{equation}
By Lemma \ref{l2.1}, 
\[
\Big|   \bens \big[ (\partial_{\!s} + nL_n) g_{s,t}(M_{ns}) \big] \Big| 
		\leq C(a,b) \Big(\frac{\| D^2 g_{s,t}\|_\infty}{n} + \frac{\| D^3 g_{s,t}\|_\infty}{n^2}+\frac{\|D^4 g_{s,t}\|_\infty}{n^2} \Big),
\]
and by Lemma \ref{l2.2}, the right-hand side of this equation is bounded by
\[
 C(a,b) \Big(\frac{\| D^2 f\|_\infty}{n} + \frac{\| D^3 f\|_\infty}{n^2}+\frac{\|D^4 f\|_\infty}{n^2} \Big)e^{-(a+b)(t-s)}.
 \]
 Therefore, the expression in \eqref{caldera} is bounded by
 \[
 C(a,b) \Big(\frac{\| D^2 f\|_\infty}{n} + \frac{\| D^3 f\|_\infty}{n^2}+\frac{\|D^4 f\|_\infty}{n^2} \Big)(1-e^{-(a+b)t}).
 \]
 Using the estimate $1-e^{-x} \leq \min\{1,x\}$, we conclude that
 \[
 \big| \bens [f(M_{tn})] - \bb E[f(x_t(M_0)] \big| \leq C(a,b) \Big(\frac{\| D^2 f\|_\infty}{n} + \frac{\| D^3 f\|_\infty}{n^2}+\frac{\|D^4 f\|_\infty}{n^2} \Big) \min\{1,t\}.
 \]
Using Lemma \ref{reglem}, the theorem is proved.
\end{proof}

\subsection{The discrete approximation}

Let $f \in \mc C_b^4(\bb R)$. Our aim is to approximate $n L_n f(M)$ by $(\Lambda_{a,b}f)(M)$. We have the following estimate:

\begin{lemma}
\label{l2.1}
There exists a finite constant $C =C(a,b)$ such that for every $n \in \bb N$ and every $f \in \mc C_b^4([0,1]; \bb R)$,
\[
\big| nL_n f(M) -\Lambda_{a,b}f(M)\big| \leq C(a,b) \Big( \frac{\|f''\|_\infty}{n} + \frac{\|f'''\|_\infty}{n^2}+ \frac{\|f^{(4)}\|_\infty}{n^2}\Big)
\]
\end{lemma}

\begin{proof}

By Taylor's formula with Lagrange rest, there exists $\theta = \theta(M)$ such that
\[
\begin{split}
n L_n f(M) 
		&= \sumxln c_x(\eta) \Big( f \big( M + \tfrac{1-2 \eta_x}{n} \big) -f(M) \Big)\\
		&= \sumxln c_x(\eta) \Big( \frac{1-2\eta_x}{n} f'(M) + \frac{1}{2 n^2} f''(M) + \frac{1-2\eta_x}{6 n^3} f'''(M) + \frac{1}{24 n^4} f^{(4)}(\theta) \Big).
\end{split}
\]
Observe now that
\[
\begin{split}
\sumxln c_x(\eta) (1-2\eta_x) 
		&= \sumxln \big( (1-\eta_x)(a+X) -\eta_x(b+n-X)\big) \\
		&= (n-X)(a+X) - X(b+n-X)\\
		&=an -(a+b) X,
\end{split}
\]
and that
\[
\begin{split}
\sumxln c_x(\eta) 
		&= \sumxln \big( (1-\eta_x)(a+X) + \eta_x(b+n-X)\big)\\
		&= (n-X)(a+X) + X(b+n-X)\\
		&= 2X (n-X) + a(n-X) +bX.
\end{split}
\]
Let us define $F,G: [0,1] \to \bb R$ as
\begin{equation}
\label{vallenar}
F(m) := a-(a+b)m, \quad G(m) := 2m(1-m) + \tfrac{1}{n} (a(1-m) +bm)
\end{equation}
for every $m \in [0,1]$. We observe that
\begin{equation}
\label{calama}
n L_n f(M) = n F(M) \Big( \frac{f'(M)}{n} +\frac{f'''(M)}{6n^3}\Big)
		+ n^2 G(M) \Big( \frac{f''(M)}{2 n^2} + \frac{f^{(4)}(M)}{24 n^4}\Big).
\end{equation}
Therefore,
\[
\big| nL_n f(M) -\Lambda_{a,b}f(M)\big| = \Big| \frac{F(M) f'''(M)}{6 n^2}+\frac{a(1-M)+bM}{2n} f''(M) + \frac{G(M)}{24 n^2} f^{(4)}(M)\Big|.
\]
Estimating the derivatives by their uniform norms, the lemma is proved.
\end{proof}

\subsection{Regularity of Fokker-Planck equation}

Our aim in this section is to prove that the solutions of the Fokker-Planck equation \eqref{fokpla} with smooth final conditions are smooth. More precisely, we will prove the following result:

\begin{lemma}
\label{l2.2}
Let $m \in \bb N_0$, let $f \in \mc C_b^m([0,1], \bb R)$, let $t >0$ and let $(g_{s,t}; s \in [0,t])$ be the solution of \eqref{fokpla} with final condition $f$. For every $\ell \in \{0,1,\dots,m\}$ and every $s \in [0,t]$,
\[
\|\partial_{\!x}^\ell g_{s,t}\|_\infty \leq \| \partial_{\!x}^\ell f\|_\infty e^{-\ell(a+b+\ell-1)(t-s)}.
\]
In particular, for every $\ell \in \{1,\dots,m\}$,
\[
\|\partial_{\!x}^\ell g_{s,t}\|_\infty \leq \| \partial_{\!x}^\ell f\|_\infty e^{-(a+b)(t-s)}.
\]
\end{lemma}

\begin{proof}
Let $f \in \mc C^\infty_b([0,1]; \bb R)$, and let $(u_t(x); t \geq 0)$ be the solution of the equation
\begin{equation}
\label{copiapo}
\partial_{\!t} u = \Lambda_{a,b} u
\end{equation}
with initial condition $u_0 = f$. Let $(x_t(m); t \geq 0)$ the solution of \eqref{wrifish}. By It\^o's formula applied to $u_{t-s}(x_s(m))$, $u$ admits the representation
\[
u_t(m) = \bb E \big[ f \big(x_t(m)\big) \big].
\]
In particular, 
\[
\|u_t\|_\infty \leq \|f\|_\infty 
\]
for every $t \geq 0$. Now observe that 
\[
\partial_{\!x} \Lambda_{a,b} = \Lambda_{a,b} \partial_{\!x} + (1-2x) \partial^2_{\!x} - (a+b) \partial_{\!x} = \big( \Lambda_{a+1,b+1} - (a+b) \big) \partial_{\!x}.
\]
Using this relation, we see by induction that for every $\ell \in \bb N$,
\[
\partial_{\!t} \big( \partial_{\!x}^\ell u\big) = \Lambda_{a+\ell,b+\ell}\big( \partial_{\!x}^\ell u\big) - \ell (a+b+\ell-1) \partial_{\!x}^\ell u.
\]
For every $\alpha, \beta >0$ and every $\lambda \in \bb R$, the solution of the equation
\[
\partial_{\!t} v = \Lambda_{\alpha,\beta} v - \lambda v
\]
admits the representation
\[
v_t(x) = \bb E \big[f(y_t(x))e^{-\lambda t} \big],
\]
where $(y_t(x); t \geq 0)$ is the solution of \eqref{wrifish} with $a =\alpha$, $b = \beta$ and initial condition $x \in [0,1]$.
In particular 
\[
\|\partial_{\!x}^\ell u_t \|_\infty \leq \|\partial_{\!x}^\ell f\|_\infty e^{-\ell(a+b+\ell-1)t}.
\]
Observing that $(g_{t-s,t}; s \in [0,t])$ is a solution of \eqref{copiapo}, the lemma is proved.
\end{proof}

\section{Thermalization and convergence of local densities}

For every $\sigma \in \On$, let $u(\sigma)$ denote the uniform law in $\Omega_{n, X(\sigma)}$. By definition, $\sigma \in \Omega_{n, X(\sigma)}$. 
Let $\vec{\mathbf 0}, \vec{\mathbf{ 1}}$, be the constant configurations with values equal to $0$, $1$, respectively, and let $\Onast:= \On \setminus\{\vec{\mathbf 0}, \vec{\mathbf 1}\}$.
Let $\sigma \in \Onast$ be fixed. In this section we will compute the times at which the laws of $\eta(t)$ with respect to $\pens$ and with respect to $\penus$ become close. Observe that for $\sigma = \vec{\mathbf 0}, \vec{\mathbf 1}$ both laws are identical. 


\begin{theorem}
\label{phase1}
Fix $\alpha <\frac{1}{3}$ and for every $n \in \bb N$ let 
\[
B_n^\alpha := \big\{ \sigma \in \Onast; M(\sigma) (1 -M(\sigma)) \geq n^{-\alpha} \big\}.
\]
For every $\tau \in \bb R$,
\[
\lim_{n \to \infty} \sup_{\sigma \in B_n^\alpha} \Big| \frac{1}{\sqrt n} d_K(\sigma P_{\!t_n}, u(\sigma) P_{\!t_n}) - 2e^{-\tau} \Big| = 0,
\]
where
\[
t_n := \frac{1}{2} \log n + \log M(\sigma)(1-M(\sigma)) + \tau.
\]
\end{theorem}

\subsection{The local densities}

For $\sigma \in \Onast$, let $(A^i \! (\sigma); i \in Q)$ be the non-trivial partition of $\Ln$ given by
\[
A^i \! (\sigma) := \{ x \in \Ln; \sigma_x =i\} \; \forall  i \in Q,
\]
and let $\vecX(\sigma): \On \to [0,\infty)^Q$ be given by
\[
X^i \! (\eta; \sigma) := \!\! \sum_{x \in A^i \! (\sigma)} \! \! \eta_x
\]
for every $i \in Q$ and every $\eta \in \On$. Let $\vq$ be the space of probability measures in $Q$. Observe that $\vq \subseteq [0,1]^Q$.
Let $\va(\sigma) \in \vq$ be given by
\[
a^i \! (\sigma) := \frac{\# A^i \! (\sigma)}{n} \; \forall i \in Q,
\]
and let $\vecM(\sigma) : \On \to [0,1]^Q$ be given by
\[
M^i \! (\eta;\sigma) := \frac{X^i \! (\eta;\sigma)}{a^i \! (\sigma) n}
\]
for every $i \in Q$ and every $\eta \in \On$. We will call the vector $\vecM$ the vector of \emph{local densities} associated to the partition generated by $\sigma$. 

Since $\sigma$ is going to be kept fixed most of the time, we will exclude the dependence on $\sigma$ from the notation, unless we need to keep explicit track of it. 

We say that a measure $\mu$ in $\On$ is \emph{$\sigma$-invariant} if 
\[
\mu(\eta) = \mu(\eta') \text{ whenever } \vecX( \eta; \sigma) =\vecX( \eta' ; \sigma).
\]
Observe that $\sigma$-invariant measures are determined by the law of $\vecX(\sigma)$.
This is the case for the law of $\eta(t)$ under $\pens$. Therefore, $\Law(\eta(t); \pens)$ is determined by $\Law(\vecX(\eta(t)); \pens)$ and equivalently by $\Law(\vecMt; \pens)$. The same is true for the law of $\eta(t)$ under $\penus$. The following result establishes that in order to compute the Kantorovich distance between $\sigma P_{\!t}$ and $u(\sigma) P_{\!t}$, it is enough to compute the distance between the corresponding local densities:

\begin{lemma}
\label{acopl}
Let $\sigma \in \Onast$ and let $\mu$, $\nu$ be $\sigma$-invariant measures in $\On$. We have that
\[
d_K(\mu,\nu) = d_K \big(\Law(\vecX ; \mu), \Law(\vecX; \nu) \big).
\]
\end{lemma}

\begin{proof}
Since $\vecX(\sigma)$ is a Lipschitz function with Lipschitz constant equal to $1$, the inequality
\[
d_K \big(\Law(\vecX ; \mu), \Law(\vecX; \nu) \big) \leq d_K(\mu,\nu)
\]
follows from general properties of the Kantorovich distance. In order to prove the converse estimate, let $\pi(\vec x, \vec y)$ be a coupling of the laws of $\vecX$ under $\mu$ and $\nu$. We construct a coupling between $\mu$ and $\nu$ as follows. Observe that $\pi$ is a discrete measure. Fix $\vec x$, $\vec y$ such that $\pi(\vec x, \vec y) >0$. If $x^i \leq y^i$, we put $y^i$ particles in $A^i$ uniformly at random, which determines the values of $\eta'$ in $A^i$. Then we choose $x^i$ of these particles uniformly at random, which determines the values of $\eta$ in $A^i$. If $x^i > y^i$, we put $x^i$ particles in $A^i$ uniformly at random, which determines the values of $\eta$ in $A^i$. Then we choose $y^i$ of these particles uniformly at random, which determines the values of $\eta'$ in $A^i$. We do this independently at each $A^i$. 
The measure in $\On \times \On$ obtained in this way is a coupling of $\mu$ and $\nu$ satisfying
\[
\bb H_n(\eta, \eta') = \| \vec x -\vec y \|
\]
for every $\vec x$, $\vec y$ with $\pi(\vec x, \vec y) \neq 0$. 
By the Kantorovich-Rubinstein duality, we conclude that
\[
d_K(\mu,\nu) \leq d_K \big(\Law(\vecX ; \mu), \Law(\vecX; \nu) \big),
\]
which proves the lemma. 
\end{proof}

Thanks to this lemma, we see that the proof of Theorem \ref{phase1} is reduced to the analysis of the behavior of the local densities with respect to the laws $\pens$  and $\penus$.

\subsection{\emph{A priori} estimates}

Our aim in this section is to prove a quantitative law of large numbers under the measure $\pens$ for the vector $\vecMt := \vecM (\eta(t); \sigma)$. This quantitative bound will be used later on as an input on the proof of the QCLT for the local densities. For this reason, we call this estimate an \emph{a priori estimate}. First we start proving an \emph{a priori} estimate for $M_t$. For every $m \in [0,1]$, let $(m_t(m); t \geq 0)$ be the solution of the equation
\begin{equation}
\label{EDO}
\tfrac{d}{dt} m_t = \tfrac{1}{n} F(m_t)
\end{equation}
with initial condition $m_0 =m$. We have the following result:

\begin{lemma}[First \emph{a priori} estimate] 
\label{l3.1}
There exists a finite constant $C = C(a,b)$ such that for every $n \in \On$, every $\sigma \in \On$ and every $t \geq 0$, 
\[
\bens \big[ \big( M_t -m_t\big)^2\big] \leq \frac{C(a,b)}{n} \min\{n,t\} \leq \frac{C(a,b)t}{n},
\]
where $(m_t; t \geq 0)$ is the solution of \eqref{EDO} with initial condition $m_0 = M(\sigma)$.
\end{lemma}

\begin{proof}
The following formula is very useful: for every $f: \On \to \bb R$,
\[
L_n f^2 = 2 f L_n f+ \Gamman f.
\]
Similarly, if $f: \On \times [0,T] \to \bb R$ is differentiable in $t \in [0,T]$, then
\[
(\partial_{\!t} + L_n) f^2 = 2 f (\partial_{\!t} + L_n) f + \Gamman f.
\]
Taking $f_1(x) =\frac{x}{n}$ and $f_2(x) = \frac{x^2}{n^2}$ in \eqref{calama}, we see that
\[
L_n M = \tfrac{1}{n} F(M)
\]
and
\[
\Gamman M = \frac{1}{n} G(M),
\]
where $F,G$ were defined in \eqref{vallenar}. Therefore,
\[
\begin{split}
(\partial_{\!t} +L_n)(M-m_t)^2 
		&= 2 (M-m_t) (\partial_{\!t} + L_n) (M-m_t) + \Gamman M \\
		&=\tfrac{2}{n}(M-m_t) \big( F(M) -F(m_t) \big) + \tfrac{1}{n}G(M)\\
		&= - \tfrac{2(a+b)}{n} (M-m_t)^2 +\tfrac{1}{n} G(M).
\end{split}
\]
By Dynkin's formula, 
\[
\tfrac{d}{dt} \bens \big[ \big( M_t - m_t\big)^2 \big] \leq - \tfrac{2(a+b)}{n} \bens \big[ \big( M_t-m_t\big)^2 \big] + \tfrac{1}{n} \bens \big[G(M_t)\big].
\]
By Gronwall's inequality, recalling that $m_0 =M_0= M(\sigma)$, we see that
\[
\bens \big[ \big( M_t - m_t\big)^2 \big] \leq \frac{\|G\|_\infty}{2(a+b)} \big( 1- e^{-( \frac{2(a+b)}{n}) t} \big).
\]
Using the bound $1-e^{-x} \leq \min\{1,x\}$, the lemma is proved. 
\end{proof}

In view of Lemma \ref{l3.1}, it is natural to define the \emph{density fluctuation process} $(\xi_t; t \geq 0)$ as
\[
\xit := \frac{1}{\sqrt n} \sum_{x \in \Ln} \big( \eta_x(t) - m_t\big)
\]
for every $t \geq 0$. The estimate of Lemma \ref{l3.1} can be rewritten as
\[
\bens[\xit^2] \leq C t
\]
for every $t \geq 0$. Observe that $\xit = \sqrt n(M_t-m_t)$.

Observe that equation \eqref{EDO} is linear, and therefore it can be solved explicitly: we have that
\[
m_t = m_0 e^{- (\frac{a+b}{n}) t} + \tfrac{a}{a+b} \big( 1- e^{- (\frac{a+b}{n}) t} \big),
\]
and in particular
\begin{equation}
\label{coquimbo}
\begin{split}
|m_t - m_0| 
		&= \big| \big( \tfrac{a}{a+b} -m_0 \big) \big( 1 - e^{- (\frac{a+b}{n}) t}\big) \big|\\
		&\leq \frac{C(a,b)}{n} \min\{n,t\} \leq \frac{C(a,b)t}{n}.
\end{split}
\end{equation}
In other words, $m_t$ does not move macroscopically for times $t \ll n$. 

Once we have proved an estimate for $\bens[M_t^2]$, we can derive an effective estimate for $\bens[ (M_t^i)^2]$. We have the following result:

\begin{lemma}[Second \emph{a priori} estimate]
\label{l3.2}
There exists a finite constant $C = C(a,b)$ such that for every $n \in \bb N$, every $\sigma \in \Onast$ and every $t \geq 0$,
\[
\bens\big[ \big(M_t^i-m_t^i\big)^2 \big] \leq \frac{C t}{a^i \! n}.
\]
\end{lemma}

\begin{proof}
Let us compute $L_n M^i$:
\[
\begin{split}
L_n M^i 
	&= \frac{1}{a^i \! n^2} \sumxai \big( (1-\eta_x)(a+X) - \eta_x (b+n-X) \big)\\
	&= \frac{1}{a^i \! n^2} \big( (a^i \! n - X^i)(a+X) - X^i \! (b+n-X) \big)\\
	&= (1-M^i) \big(\tfrac{a}{n} + M \big) - M^i \big( \tfrac{b}{n} +1 -M \big)\\
	&= \tfrac{a}{n} + M - \big(1+\tfrac{a+b}{n}\big) M^i.
\end{split}
\]
Observe that $M = \va \cdot \vecM$. Let us define $\vecF(\va): [0,1]^Q \to \bb R^Q$ as 
\[
F^i \! (\vecm; \va) := \tfrac{a}{n} + \va \cdot \vecm - \big(1+ \tfrac{a+b}{n} \big) m^i
\]
for every $i \in Q$ and every $ \vecm \in [0,1]^Q$. We have that $L_n M^i = F^i \! (\vecM; \va)$. Taking in consideration the proof of Lemma \ref{l3.1}, let us define $(\vecmt(\va); t \geq 0)$ as the solution of the equation
\begin{equation}
\label{EDO2}
\tfrac{d}{dt} \vecmt = \vecF(\vecmt)
\end{equation} 
with initial condition $\vecmzero := (0,1)$.

Now let us compute $\Gamman M^i$:
\[
\begin{split}
\Gamman M^i 
		&= \frac{1}{(a^i)^2n^3} \sumxai \big( (1-\eta_x)(a+X) + \eta_x(b+n-X)\big) \\
		&= \frac{1}{a^i \! n} \big( (1-M^i) \big( \tfrac{a}{n} +M \big) + M^i \big( \tfrac{b}{n} +1 -M \big)\\
		&= \frac{1}{a^i \! n} \big( M+M^i- 2 MM^i + \tfrac{1}{n} ( a(1-M^i) + b M^i ) \big).
\end{split}
\]
Let us define $\vecG(\cdot; \va) : [0,1]^Q \to \bb R^Q$ as
\[
\label{defgi}
G^i \! (\vecm; \va) := \va \cdot \vecm + m^i -2(\va \cdot \vecm) m^i +\tfrac{1}{n}(a(1-m^i) + bm^i)
\]
for every $i \in Q$ and every $\vecm \in [0,1]^Q$. We see that $\Gamman M^i = \frac{1}{a^i \! n} G^i \! (\vec M)$. Therefore,
\[
\begin{split}
(\partial_{\!t}+L_n) (M^i-m_t^i)^2 
		&= 2(M^i-m_t^i) \big( F^i \! (\vecM) - F^i \! (\vecmt) \big) + \frac{G^i \! (\vecM)}{a^i \! n} \\
		&= -2 \big(1+\tfrac{a+b}{n} \big) (M^i-m_t^i)^2 + 2(M^i-m_t^i) ( M -m_t) + \frac{G^i \! (\vecM)}{a^i \! n}.
\end{split}
\]
By the weighted Cauchy-Schwartz inequality, for every $\delta >0$ we have the estimate
\begin{equation}
\label{losandes}
(\partial_{\!t}+L_n) (M^i-m_t^i)^2 \leq -2\big(1+\tfrac{a+b}{n}-\delta\big) (M^i-m_t^i)^2 + \frac{(M-m_t)^2}{2 \delta} + \frac{\|G^i\|_\infty}{a^i \! n}.
\end{equation}
By Dynkin's formula and Gronwall's inequality, we see that for every $\delta \in (0,1]$,
\[
\bens \big[ \big(M_t^i-m_t^i\big)^2 \big] \leq \frac{1}{2\delta}\int_0^t e^{-2(1+\tfrac{a+b}{n}-\delta)(t-s)} \bens\big[\big(M_s-m_s\big)^2 \big] ds
		+ \frac{\|G^i\|_\infty}{a^i \! n} \frac{1-e^{-2(1+\tfrac{a+b}{n} -\delta)t}}{2\big(1+ \tfrac{a+b}{n} -\delta \big)}.
\]
Using the estimate of Lemma \ref{l3.1}, we obtain the estimate
\[
\bens \big[ \big(M_t^i-m_t^i\big)^2 \big]  \leq \Big( \frac{C t}{\delta n} + \frac{\|G^i\|_\infty}{a^i \! n} \Big) \min\Big\{ t, \frac{1}{2 \big(1+\tfrac{a+b}{n}-\delta \big) } \Big\}.
\]
Observing that there exists a finite constant $C = C(a,b)$ such that $\|G^i \! (\va)\|_\infty \leq C$ for every $i \in Q$ and every $\va \in \vq$, the lemma is proved.
\end{proof}

In view of Lemma \ref{l3.2}, it is natural to introduce the \emph{density fluctuation process} as follows. Let $\vecxi: \On \times \Onast \times [0,\infty)$ be given by
\[
\xi^i \! (\eta; \sigma, t) := \sqrt n \Big( \frac{X^i \! (\eta; \sigma)}{n} - a^i \! (\sigma)  m_t^i(\va(\sigma))\Big)
\]
for every $i \in Q$, every $\eta \in \On$, every $\sigma \in \Onast$ and every $t \geq 0$. Then we define the process $(\vecxit; t \geq 0)$ as $\xit^i := \xi^i \! (\eta(t); \sigma,t)$ for every $t \geq 0$ and every $i \in Q$. Observe that the estimate of Lemma \ref{l3.2} can be rewritten as
\[
\bens \big[ \big(\xit^i\big)^2 \big] \leq C t
\]
for every $i \in Q$ and every $t \geq 0$. Observe that
\[
\xit^i = a^i \sqrt{n} \big( M_t^i -m_t^i\big)
\]
and that
$
\xit = \vec{1} \cdot \vecxit.
$
\subsection{The discrete approximation} Let $f \in \mc C_b^4(\bb R^Q; \bb R)$. Our aim in this section is to show that $(\partial_{\!t} +L_n) f(\vecxit)$ can be well approximated by a linear operator applied to $f$ and evaluated at $\vecxit$. For $f \in \mc C_b^2(\bb R^Q; \bb R)$ and $t \geq 0$, let $\bb L_t f: \bb R^Q \to \bb R$ be given by
\[
\bb L_t f(\vecx) := \sumiq \Big( \big(a^i(\vec{1} \cdot \vecx) -\big( 1+\tfrac{a+b}{n} \big) x^i \big) \partial_{\hspace{-0.6pt}i} f(\vecx) + \frac{a^i G^i \! (\vecmt)}{2} \partial_{\hspace{-0.6pt}i}^2 f(\vecx) \Big)
\]
for every $\vecx \in \bb R^Q$. We have the following estimate:

\begin{lemma}[Discrete approximation]
\label{l3.3} 
There exists a finite constant $C = C(a,b)$ such that
\[
\big| (\partial_{\!t} +L_n) f(\vecxit) - \bb L_t f( \vecxit) \big| \leq C \sumiq \Big( \frac{\|\partial_{\!i}^2 f\|_\infty\|\vecxit\|_{\ell^1}}{\sqrt n} + \frac{\|\partial_{\!i}^3 f\|_\infty}{\sqrt n} + \frac{\|\partial_{\!i}^4 f\|_\infty}{n} \Big)
\]
for every $n \in \bb N$, every $t \geq 0$, every $\sigma \in \Onast$ and every $f \in \mc C_b^4(\bb R^Q; \bb R)$.
\end{lemma}

\begin{proof}
Let us compute $L_nf(\vecxit)$. By Taylor's formula of order $3$ with Lagrange rest, there exist vectors $\vec{\theta}_{\!t,x}^i$ such that
\[
\begin{split}
L_n f(\vecxit) 
		&= \frac{1}{n} \sumiq \sumxai c_x(\eta) \Big( f\big( \vecxit + \tfrac{(1-2\eta_x)e_i}{\sqrt n} \big) - f ( \vecxit) \Big)\\
		&=\frac{1}{n} \sumiq \sumxai \Big[ c_x(\eta) (1-2\eta_x) \Big( \frac{1}{\sqrt n} \partial_{\!i} f(\vecxit) + \frac{1}{6 n^{3/2}} \partial_{\!i}^3 f(\vecxit) \Big)\\
		&\quad \quad \quad \quad \quad \quad +c_x(\eta) \Big( \frac{1}{2n} \partial_{\!i}^2 f(\vecxit) + \frac{1}{24 n^2} \partial_{\!i}^4 f(\vec{\theta}_{\!t,x}^i)\Big)\Big].
\end{split}
\]
Now we observe that
\[
\frac{1}{n} \sumxai c_x(\eta) (1-2\eta_x) = a^i \! n F^i \! (\vecM)
\]
and that
\[
\frac{1}{n} \sumxai c_x(\eta) = a^i \! n G^i \! (\vecM).
\]
Therefore,
\[
\begin{split}
L_n f( \vecxit) 
		&= \sumiq \Big( a^i \! n F^i \! (\vecM) \Big(  \frac{1}{\sqrt n} \partial_{\!i} f(\vecxit) + \frac{1}{6 n^{3/2}} \partial_{\!i}^3 f(\vecxit) \Big)\\
		&\quad \quad \quad \quad +\frac{a^i G^i \! (\vecM)}{2} \partial_{\!i}^2 f(\vecxit) + \frac{\mc R_{\!t}^i \! (f)}{n} \Big),
\end{split}
\]
where $\mc R_{\!t}^i \! (f)\|_\infty \leq C(a,b) a^i \|\partial_{\!i}^4 f\|_\infty$. Observe as well that
\[
\partial_{\!t} f(\vecxit) = -a^i \sqrt n \sumiq \partial_{\!i} f(\vecxit) F^i \! (\vecmt),
\]
from where
\begin{equation}
\label{olmue}
\begin{split}
\Big| (\partial_{\!t} + L_n) f(\vecxit) - \sumiq a^i \Big( \sqrt{n} \big( &
		F^i \! (\vecM) -F^i \! (\vecmt) \big) \partial_{\!i} f(\vecxit)\\
		& +\frac{G^i \! (\vecM)}{2} \partial_{\!i}^2 f(\vecxit)
 + \frac{F^i \! (\vecM)}{6 \sqrt n} \partial_{\!i}^3 f(\vecxit) \Big) \Big| 
		\leq \frac{C(a,b)}{n} \sumiq a^i \|\partial_{\!i}^4 f\|_\infty.
\end{split}
\end{equation}
Observe that there exists a finite constant $C = C(a,b)$ such that $\|F^i\|_\infty \leq C$ for every $i \in Q$ and every $\va \in \vq$. Therefore,
\begin{equation}
\label{quilpue}
\begin{split}
\Big| (\partial_{\!t} + L_n) f(\vecxit) - \sumiq a^i \Big( \sqrt{n} \big( 
		F^i \! (\vecM) -
		&F^i \! (\vecmt) \big) \partial_{\!i} f(\vecxit)
		+\frac{G^i \! (\vecM)}{2} \partial_{\!i}^2 f(\vecxit) \Big)\Big| \leq \\
		&\leq C(a,b) \sumiq a^i \Big(\frac{\|\partial_{\!i}^3 f\|_\infty}{\sqrt n} + \frac{\|\partial_{\!i}^4 f\|_\infty}{n}\Big).
\end{split}
\end{equation}
In order to rewrite this estimate in terms of the operator $\bb L_t$, observe that 
\begin{equation}
\label{algarrobo}
\begin{split}
\sqrt n \big( F^i \! (\vecM) - F^i \! (\vecmt) \big)
		&= \sqrt n (M-m_t) - \sqrt n \big( 1+\tfrac{a+b}{n} \big) (M^i -m_t^i) \\
		&= \xit - \big( 1+ \tfrac{a+b}{n} \big) \frac{\xit^i}{a^i}
\end{split}
\end{equation}
and observe that
\begin{equation}
\label{elquisco}
|G^i \! (\vecM) - G^i \! (\vecmt) | 
		\leq \sum_{j \in Q} \|\partial_{\!j} G^i\|_\infty |M^i -m_t^i| 
		\leq \sum_{j \in Q} \frac{\|\partial_{\!j} G^i\|_\infty |\xit^i|}{a^i \sqrt n}.
\end{equation}
We have that $\|\partial_{\!j} G^i\|_\infty \leq a^j +C(a,b) \delta_{ij}$, from where 
\begin{equation}
\label{eltabo}
a^i |G^i \! (\vecM) - G^i \! (\vecmt)| \leq C(a,b) \sum_{j \in Q} \frac{|\xit^j|}{\sqrt n} = \frac{C(a,b) \|\vecxit\|_{\ell^1}}{\sqrt n}. 
\end{equation}
Putting identity \eqref{algarrobo} and estimates \eqref{elquisco}, \eqref{eltabo} into \eqref{quilpue}, the lemma is proved.
\end{proof}

\subsection{Regularity of the Fokker-Planck equation}

In this section we will prove a regularity estimate for solutions of an auxiliary PDE. This estimate will be needed in order to prove a QCLT for the process $(\vecxit; t \geq 0)$. Let $f \in \mc C_b(\bb R^Q;  \bb R)$ and let $t \geq 0$. Let $(g_{s,t}; s \in [0,t])$ be the solution of the \emph{backwards Fokker-Planck equation} 
\begin{equation}
\label{fokpla2}
(\partial_{\!s} + \bb L_s) g_s = 0
\end{equation}
with final condition $g_{t,t} = f$. We have the following result. 

\begin{lemma}[Regularity estimate] 
\label{l3.4}
Let $\ell \in \bb N_0$. There exists $C = C(\ell)$ such that for every $f \in \mc C_b^\ell(\bb R^Q; \bb R)$ and every $s,t \in [0,\infty)$ such that $s \leq t$, the solution $(g_{s,t}; s \in [0,t])$ of \eqref{fokpla2} satisfies
\[
\|\partial^k_i g_{s,t} \|_\infty \leq C(\ell) \| D^k f\|_\infty 
\]
for every $i \in Q$ and every $k \in \{0,1,\dots,\ell\}$.
\end{lemma}

\begin{proof}
Let $(\vecyt; t \geq 0)$ the diffusion generated by the operators $(\bb L_t; t \geq 0)$ with initial condition $\vecyzero = \vec 0$, that is, the solution of the equations
\begin{equation}
\label{sde}
d y_{\!t}^i = \big( a^i (\vec 1 \cdot \vecyt) - \big(1+\tfrac{a+b}{n} \big) y_{\!t}^i \big) dt + \sqrt{a^i G^i \! (\vecmt)} d B_t^i, \quad \forall i \in Q,
\end{equation}
with zero initial condition, where $(\vec{B}_t; t \geq 0)$ is a standard Brownian motion in $\bb R^Q$. Applying It\^o's formula to $g_{s,t}(\vecys)$, we see that $g_{s,t}$ satisfies 
\[
g_{s,t}(\vecx) = \bb E \big[ f(\vecyt) \big| \vecys = \vecx \big],
\]
and in particular $\|g_{s,t} \|_\infty \leq \|f\|_\infty$. This proves the lemma for $\ell =0$. The idea is to compute the equation satisfied by the derivative $D g_{s,t}$ and to argue in a similar way. It will be useful to introduce the \emph{commutator} 
\[
[\partial_{\!i}, \bb L_t] := \partial_{\!i} \bb L_t - \bb L_t \partial_{\!i}
\]
between $\partial_{\!i}$ and $\bb L_t$. We have that
\[
\begin{split}
[ \partial_{\!i}, \bb L_t] 
		&= \big[ \partial_{\!i} , \sum_{j \in Q} \big( a^i (\vec 1 \cdot \vec x) - \big( 1+\tfrac{a+b}{n} \big) x^j\big)\partial_{\!j}\big]\\
		&= \sum_{j \in Q} \big( a^i - \delta_{ij} \big(1+\tfrac{a+b}{n}\big)\big) \partial_{\!j}\\
		&= a^i \sum_{j \in Q} \partial_{\!j} - \big(1+\tfrac{a+b}{n} \big) \partial_{\!i}.
\end{split}
\]
In view of this formula, let us introduce the \emph{divergence operator}
\[
\Div := \sumiq \partial_{\!i},
\]
in such a way that
\[
[\partial_{\!i}, \bb L_t] = a^i \Div - \big(1+\tfrac{a+b}{n} \big) \partial_{\!i}.
\]
This definition allows us to diagonalize the commutator $[D, \bb L_t]$. In fact,
\[
[\Div, \bb L_t] = - \big(\tfrac{a+b}{n} \big) \Div \quad \text{ and } \quad [\partial_{\!i} -a^i \Div, \bb L_t] = - \big(1+\tfrac{a+b}{n}\big) (\partial_{\!i} -a^i \Div).
\]
Define $\mc A_i := \partial_{\!i} -a^i \Div$. Observe that if $(h_s; s \in [0,t])$ satisfies
\[
(\partial_{\!s} + \bb L_s - \gamma ) h_s =0,
\]
for some $\gamma \in \bb R$, then $\mc A_i h$ and $\Div h$ satisfy
\[
\big(\partial_{\!s} + \bb L_s - \gamma -\tfrac{a+b}{n} \big) \Div h_s = 0 \quad \text{and} \quad \big(\partial_{\!s} + \bb L_s - \gamma -1-\tfrac{a+b}{n} \big) \mc A_i h_s=0.
\]
Recursively, we see that for $\ell, m \in \bb N_0$, the function $h:= \Div^m \mc A_i^\ell  g$ satisfies the equation
\[
\big( \partial_{\!s} + \bb L_s -\ell -(\ell+m) \big( \tfrac{a+b}{n} \big) \big) h =0
\]
with final condition $h_t = \Div^m \mc A_i^\ell f$. By It\^o's formula, $h$ admits the representation
\[
h_s(\vecx) = \bb E\big[ \Div^m \mc A_i^\ell f(\vecyt) e^{-\lambda_{\ell,m}(t-s)} \big| \vecys = \vecx \big],
\]
where $\lambda_{\ell,m} := \ell + (\ell+m)\frac{a+b}{n}$. Observe now that
\[
\partial_{\!i}^k = (\mc A_i + a^i \Div)^k = \sum_{j=0}^k \binom{k}{j} (a^i)^{k-j}  \Div^{k-j} \mc A_i^j,
\]
from where
\[
\begin{split}
\partial_{\!i}^k g_{s,t} (\vec x) 
		&=  \sum_{j=0}^k \binom{k}{j} (a^i)^{k-j}  \Div^{k-j} \mc A_i^j g_{s,t}(\vecx)\\
		&= \bb E\Big[  \sum_{j=0}^k \binom{k}{j} e^{-\lambda_{j,k-j}(t-s)}(a^i)^{k-j}  \Div^{k-j} \mc A_i^jf(\vecyt) \Big| \vecys = \vecx \Big] \\
		&= \bb E \big[ \big(\mc A_i e^{-(t-s)} +a^i \Div \big)^k f(\vecyt) e^{-k(\frac{a+b}{n})(t-s)} \big| \vecys = \vecx \big].
\end{split}
\]
We conclude that
\[
\big\| \partial_{\!i}^k g_{s,t} \big\|_\infty \leq \big\| \big(\mc A_i e^{-(t-s)} +a^i \Div \big)^k f\|_\infty e^{-k\frac{a+b}{n} (t-s)},
\]
which proves the lemma.
\end{proof}

\subsection{The quantitative CLT}
In this section we will prove a quantitative CLT for the fluctuation vectors $(\vecxit; t \geq 0)$. Recall the definition \eqref{sde} of $(\vecyt; t \geq 0)$. We will prove the following lemma:

\begin{lemma}
\label{l3.5}
There exists a finite constant $C(a,b)$ such that for every $n \in \bb N$, every $\sigma \in \Onast$ and every $t \leq \sqrt n$,
\[
d_K(\vecxit,  \vecyt) \leq \frac{C(t^{1/4}+t^{3/4})}{n^{1/6}}.
\]
\end{lemma}

\begin{proof}
The proof of this lemma follows a parabolic version of \emph{Stein's method}. Recall that we want estimate the difference
\[
\bb E_n^\sigma \big[f \big(\vecxit \big)\big] - \bb E[f(\vecyt)]
\]
in terms of the norm $\|Df\|_\infty$. Thanks to the Regularization Lemma \ref{reglem}, it is enough to estimate this difference in terms of norms of derivatives of $f$ of higher order. Let $(g_{s,t}; s \in [0,t])$ be the solution of \eqref{fokpla2} with final condition $f$. Using It\^o's formula for $g_{s,t}(\vecys)$, we see that
\begin{equation}
\label{casablanca}
\bb E[f(\vecyt)] = g_{0,t}(\vec 0).
\end{equation}
Using Dynkin's formula for $g_{s,t}(\vecxis)$, we see that
\begin{equation}
\label{lascruces}
\bb E_n^\sigma[ f(\vecxit)] = g_{0,t}(\vec 0) + \int_0^t \bb E_n^\sigma \big[ (\partial_{\!s} +L_n) g_{s,t}(\vecxis) \big] ds
\end{equation}
Observe that the partial derivative $\partial_{\!s}$ acts in both $g_{s,t}$ and $\vecxis$. Therefore,
\[
(\partial_{\!s} + L_n) g_{s,t}(\vecxis) = (\partial_{\!u} + L_n) g_{s,t}(\vecxiu) \big|_{u=s} - \bb L_s g_{s,t}(\vecxis).
\]
By Lemma \ref{l3.3}, we see that
\[
|(\partial_{\!s} +L_n) g_{s,t}(\vecxis) |
		\leq C(a,b) \sumiq \Big( \frac{\|\partial_{\!i}^2 g_{s,t}\|_\infty \|\vecxis\|_{\ell^1}}{\sqrt n} + \frac{\|\partial_{\!i}^3 g_{s,t}\|_\infty}{\sqrt n} + \frac{\|\partial_{\!i}^4 g_{s,t}\|_\infty}{n}\Big).
\]
By Lemma \ref{l3.4}, we see that
\begin{equation}
\label{cartagena}
|(\partial_{\!s} +L_n) g_{s,t}(\vecxis) |
		\leq C(a,b) \sumiq \Big( \frac{\|\partial_{\!i}^2 f\|_\infty \|\vecxis\|_{\ell^1}}{\sqrt n} + \frac{\|\partial_{\!i}^3 f\|_\infty}{\sqrt n} + \frac{\|\partial_{\!i}^4 f\|_\infty}{n}\Big).
\end{equation}
By Lemma \ref{l3.2}, 
\begin{equation}
\label{sanantonio}
\bb E_n^\sigma [ \|\vecxis\|_{\ell^1}] \leq \sumiq \bb E_n^\sigma [(\xi_{\!s}^i)^2]^{1/2} \leq C(a,b) \sumiq \sqrt{ a^i s} \leq C(a,b) \sqrt s.
\end{equation}
Using estimates \eqref{sanantonio} and \eqref{cartagena} into \eqref{lascruces} and \eqref{casablanca}, we see that
\[
\begin{split}
\big| \bb E_n^\sigma [f(\vecxit)] - \bb E[ f(\vecyt)] \big| 
		&\leq \int_0^t C(a,b)  \sumiq \Big( \frac{\|\partial_{\!i}^2 f\|_\infty \sqrt s}{\sqrt n} + \frac{\|\partial_{\!i}^3 f\|_\infty}{\sqrt n} + \frac{\|\partial_{\!i}^4 f\|_\infty}{n}\Big) ds\\
		&\leq C(a,b)  \sumiq \Big( \frac{\|\partial_{\!i}^2 f\|_\infty t^{3/2}}{\sqrt n} + \frac{\|\partial_{\!i}^3 f\|_\infty t}{\sqrt n} + \frac{\|\partial_{\!i}^4 f\|_\infty t}{n}\Big)\\
		&\leq  C(a,b)  \Big( \frac{\|D^2 f\|_\infty t^{3/2}}{\sqrt n} + \frac{\|D^3 f\|_\infty t}{\sqrt n} + \frac{\|D^4 f\|_\infty t}{n}\Big).
\end{split}
\]
By the Regularization Lemma \ref{reglem}, we conclude that
\[
d_K(\vecxit, \vecyt) \leq C(a,b) \Big( \frac{t^{3/4}}{n^{1/4}}+ \frac{t^{1/3}}{n^{1/6}} + \frac{t^{1/4}}{n^{1/4}}\Big).
\]
In this estimate, the worst term in $n$ is $n^{1/6}$. For $t \leq 1$, the worst term in $t$ is $t^{1/4}$ and for $t \geq 1$, the worst term in $t$ is $t^{3/4}$. We conclude that
\[
d_K(\vecxit, \vecyt) \leq \frac{C(a,b) (t^{1/4}+t^{3/4})}{n^{1/6}},
\]
which proves the lemma.
\end{proof}

\subsection{The uniform initial condition}

Once we have proved a quantitative CLT for the local density with respect to $\bb P_n^\sigma$, our aim is to prove a quantitative CLT for the local densities, but now under $\bb P_n^{u(\sigma)}$. In order to do that, first we observe that
\[
\bb E_n^{u(\sigma)} [\eta_x(0) ] = M(\sigma).
\]
Let $(m_t; t \geq 0)$ be the solution of \eqref{EDO} with initial condition $m_0 = M(\sigma)$ and let $(\veczetat; t \geq 0)$ be given by
\[
\zetat^i:= \frac{1}{\sqrt{n}} \big( X_t^i - a^i \! n m_t^i\big)
\]
for every $t \geq 0$ and every $i \in Q$. The plan is to repeat for $(\veczetat; t \geq 0)$ the steps leading to Lemma \ref{l3.5}. First we observe that $\veczetazero \neq 0$. Therefore, in order to use the proof of Lemma \ref{l3.2} to obtain an \emph{a priori} estimate for $\veczetat$, first we need to compute $\bb E_n^{u(\sigma)}\big[ \big( \zeta_{\!0}^i \big) ^2 \big]$ for $i \in Q$:

\begin{lemma}
\label{l3.5.1} For every $n \in \bb N$, every $\sigma \in \Onast$ and every $i \in Q$,
\[
\bb E_n^{u(\sigma)} \big[ \big( \zeta_{\!0}^i \big) ^2 \big] = \frac{n}{n-1} m_0^2(1-m_0)^2.
\]
\end{lemma}

\begin{proof}
 Observe that for $i=1$, since $\bb E_n^{u(\sigma)} [\zeta_{\!0}^i]=0$ and $ m_0 = a^1 = \frac{X(\sigma)}{n}$,
\begin{equation}
\label{angostura}
\begin{split}
\bb E_n^{u(\sigma)} \big[ \big( \zeta_{\!0}^i \big) ^2 \big]  
		&= \frac{1}{n} \sum_{\substack{x, y\in A^i \\ y \neq x}} \bb E_n^{u(\sigma)}[\eta_x \eta_y] + \frac{1}{n} \sum_{x\in A^i} \bb E_n^{u(\sigma)} [\eta_x] - \frac{X(\sigma)^4}{n^3}\\
		&= \frac{X(\sigma)(X(\sigma)-1)}{n} \cdot \frac{X(\sigma)(X(\sigma)-1)}{n(n-1)} + \frac{X(\sigma)^2}{n^2} - \frac{X(\sigma)^4}{n^3}\\
		&= \frac{X(\sigma)^2(n-X(\sigma))^2}{n^3(n-1)} = \frac{n}{n-1} m_0^2(1-m_0)^2.
\end{split}
\end{equation}
Since $\zeta_0^0 = -\zeta_0^1$, the lemma is proved.
\end{proof}

Now we have the elements needed to prove the \emph{a priori} estimate:


\begin{lemma}[Second \emph{a priori} estimate, v2] 
\label{l3.6}
There exists a finite constant $C(a,b)$ such that for every $n \in \bb N$, every $\sigma \in \Onast$, every $i \in Q$ and every $t \geq 0$,
\[
\bb E_n^{u(\sigma)}\big[ \big( \zetat^i \big) ^2 \big]  \leq C(a,b) a^i (1+t).
\]
\end{lemma}

\begin{proof}
First we observe that estimate \eqref{losandes} holds true regardless of the value of $\vecm_0$. Therefore,
\[
(\partial_{\!t}+L_n)(\zetat^i)^2 \leq -2\big( 1+ \tfrac{a+b}{n} -\delta \big) (\zetat^i)^2 + \frac{(a^i)^2\xit^2}{2 \delta} + a^i \|G^i\|_\infty.
\]
Using Dynkin's formula and Gronwall's inequality, we conclude that
\[
\begin{split}
\bb E_n^{u(\sigma)}\big[ \big( \zetat^i\big)^2\big]
		&\leq \bb E_n^{u(\sigma)}\big[ \big( \zeta_{\!0}^i\big)^2\big] e^{-2(1+\frac{a+b}{n}-\delta)t} + 
			\frac{(a^i)^2}{2 \delta} \int_0^t e^{-2(1+\frac{a+b}{n}-\delta)(t-s)} \bb E_n^{u(\sigma)} \big[\zeta_{\!s}^2 \big] ds \\
		&\quad \quad + a^i \|G^i\|_\infty  \frac{1-e^{-2(1+\tfrac{a+b}{n} -\delta)t}}{2\big(1+ \tfrac{a+b}{n} -\delta \big)}\\
		&\leq \bb E_n^{u(\sigma)}\big[ \big( \zeta_0^i\big)^2\big] + C(a,b) a^i \! t.
\end{split}
\]
Using Lemma \ref{l3.5.1} and observing that $m_0^2(1-m_0)^2 \leq a^i$ for every $i \in Q$, the lemma is proved.
\end{proof}

Now we observe that Lemma \ref{l3.3} also holds for $\veczetat$, if we replace in the definition of $\bb L_t$ the vector $\vecmt$ by the vector $(m_t,m_t)$. More precisely, let $\tilde{\bb L}_t$ be defined as
\[
\tilde{\bb L}_t f(\vecx) := \sumiq\Big( a^i (\vec{1} \cdot \vecx) - \big( 1+ \tfrac{a+b}{n} \big) x^i \big) \partial_{\!i} f(\vecx) + \frac{a^i G^i \! (m_t,m_t)}{2} \partial_{\!i}^2 f(\vecx)
\]
for every $f \in \mc C^2_b(\bb R^Q; \bb R)$ and every $\vecx \in \bb R^Q$. We have the following results:

\begin{lemma}
\label{l3.6.1} 
There exists a finite constant $C(a,b)$ such that
\[
\big| (\partial_{\!t} + L_n) f(\veczetat) - \tilde{\bb L}_t f(\veczetat) \big| \leq C(a,b) \sumiq \Big( \frac{\|\partial_{\!i}^2 f\|_\infty \|\veczetat\|_{\ell^1}}{\sqrt n} + \frac{\|\partial_{\!i}^3 f\|_\infty}{\sqrt n} + \frac{\|\partial_{\!i}^4 f\|_\infty}{n} \Big)
\]
for every $n \in \bb N$, every $\sigma \in \Onast$, every $f \in \mc C^4_b(\bb R^Q; \bb R)$ and every $t \geq 0$.
\end{lemma}

\begin{lemma}
\label{l3.6.11} Let $\ell \in \bb N_0$. There exists a finite constant $C(\ell)$ such that for every $f \in \mc C^\ell_b(\bb R^Q; \bb R)$ and every $s,t \in [0,\infty)$ such that $s \leq t$, the solution of the backwards Fokker-Planck equation
\begin{equation}
\label{fokpla3}
(\partial_{\!s} + \tilde{\bb L}_s) g_s = 0
\end{equation}
with final condition $g_{t,t}=f$ satisfies
\[
\|\partial_{\!i}^k g_{s,t} \|_\infty \leq C(\ell) \| D^k f\|_\infty
\]
for every $i \in Q$ and every $k \in \{0,1,\dots,\ell\}$.
\end{lemma}

Since the proofs of these lemmas are exactly the same of Lemmas \ref{l3.3} and \ref{l3.4}, we omit them.

Before entering into the proof of a version of Lemma \ref{l3.5} for $\veczetat$, we need to prove a quantitative CLT for $\veczetazero$. This is the content of the next section.

\subsection{The quantitative CLT for the uniform initial condition}
The aim in this section is to prove a quantitative CLT for $\veczetazero$ with respect to $\bb P_n^{u(\sigma)}$. By the particle-hole symmetry of the uniform measure, we can assume $X(\sigma) \leq \frac{n}{2}$. Since
\[
\sumiq \zeta^i_{\!0} =0,
\]
it is enough to prove a quantitative CLT for $\zeta_0^1$. In order to simplify the notation, let us write $Y:=X^1$, $\ell := X(\sigma)$ and $Z := \zeta_0^1$. Observe that $m_0 = \frac{\ell}{n}$ and that
\[
Z = \frac{1}{\sqrt n} \Big( Y-\frac{\ell^2}{n} \Big) = \frac{1}{\sqrt n} \big( Y - \ell m_0 \big).
\]
The idea is the following. We will define a Markov generator for which $u(\sigma)$ is invariant, and we will show that the action of such operator over functions of $Y$ can be well approximated by a differential operator. 

For $x,y \in \Lambda_n$ and $\eta \in \Omega_{n,\ell}$, let $\eta^{x,y} \in \Omega_{n,\ell}$ be given by
\[
\eta_z^{x,y} :=
\left\{
\begin{array}{r@{\;;\;}l}
\eta_y & z=x,\\
\eta_x & z =y,\\
\eta_z & z \neq x,y.
\end{array}
\right.
\]
For $f: \Omega_{n,\ell} \to \bb R$, let $L_{n,\ell}^{\mathrm{ex}} f: \Omega_{n,\ell} \to \bb R$ be given by
\[
L_{n,\ell}^{\mathrm{ex}} f (\eta) := \frac{1}{n} \sum_{x,y \in \Lambda_n} \big( f(\eta^{x,y})- f(\eta) \big) 
\]
for every $\eta \in \Omega_{n,\ell}$. The operator $L_{n,\ell}^{\mathrm{ex}}$ is known in the literature as the generator of the \emph{symmetric exclusion process} on the complete graph $\Lambda_n$. The only property of this generator that we will use here, is that
\[
\int L_{n,\ell}^{\mathrm{ex}} f d u(\sigma) =0
\]
for every $f: \Omega_{n,\ell} \to \bb R$. Let $f \in \mc C^3_b(\bb R)$ and let us compute $L_{n,\ell}^{\mathrm{ex}} f$:
\[
\begin{split}
L_{n,\ell}^{\mathrm{ex}} f(Z) 
		&= \frac{1}{n} \sum_{\substack{x \in A^1\\y \notin A^1}} \Big[ \eta_x(1-\eta_y) \big( f \big( Z - \tfrac{1}{\sqrt n}\big)- f(Z) \big) + \eta_y(1-\eta_x) \big( f\big( Z +\tfrac{1}{\sqrt n}\big) -f(Z)\big) \Big]\\
		&= \frac{Y(n-2\ell+Y)}{n} \big( f \big( Z - \tfrac{1}{\sqrt n}\big)- f(Z) \big) 
				+ \frac{(\ell-Y)^2}{n} \big( f \big( Z + \tfrac{1}{\sqrt n}\big)- f(Z) \big).
\end{split}
\]
Observe that
\[
 \frac{(\ell-Y)^2}{n} -  \frac{Y(n-2\ell+Y)}{n} = \frac{\ell^2 -n Y}{n} = - \sqrt n Z
\]
and that
\[
\frac{(\ell-Y)^2}{n} + \frac{Y(n-2\ell+Y)}{n} = n \gamma\big( \tfrac{Y}{n} \big),
\]
where
\[
\gamma(m) := \big( \tfrac{\ell}{n} \big)^2 + m \big(1-\tfrac{4\ell}{n} \big) +2 m^2.
\]
Using Taylor's formula of order two with Lagrange rest, there exists $\theta$ such that
\[
\begin{split}
L_{n,\ell}^{\mathrm{ex}} f(Z) 
		&= \Big(  \frac{(\ell-Y)^2}{n} -  \frac{Y(n-2\ell+Y)}{n} \Big) \frac{f'(Z)}{\sqrt n}  + \Big( \frac{(\ell-Y)^2}{n} + \frac{Y(n-2\ell+Y)}{n} \Big) \frac{f''(Z)}{2n}\\
		&\quad \quad + \frac{(\ell-Y)^2}{6n^{5/2}} f'''(Z+\theta) -  \frac{Y(n-2\ell+Y)}{6n^{5/2}} f'''(Z-\theta)\\
		&= -Z f'(Z) + \gamma\big( \tfrac{Y}{n} \big) f''(Z) + \frac{(\ell-Y)^2}{6n^{5/2}} f'''(Z+\theta) -  \frac{Y(n-2\ell+Y)}{6n^{5/2}} f'''(Z-\theta).\\
\end{split}
\]
Therefore,
\[
\Big| L_{n,\ell}^{\mathrm{ex}} f(Z)  - \frac{1}{2}\gamma\big(\tfrac{Y}{n}\big) f''(Z) + Z f' (Z) \Big| \leq \frac{\ell \|f'''\|_\infty}{3  n^{3/2}} = \frac{m_0 \|f'''\|_\infty}{3 \sqrt n}.
\]
By the mean value theorem, there exists $\tilde \theta \in \bb R$ such that
\[
\gamma\big( \tfrac{Y}{n} \big) = \gamma \big( \tfrac{\ell^2}{n^2} \big) + \frac{Z}{\sqrt n} \gamma'(\tilde \theta).
\]
Since $\gamma'(m) \leq 1$ for every $m \in [0,\frac{\ell}{n}]$,
\begin{equation}
\label{rancagua}
\Big| L_{n,\ell}^{\mathrm{ex}} f(Z)  - \frac{1}{2}\gamma\big(\tfrac{\ell^2}{n^2}\big) f''(Z) + Z f' (Z) \Big| 
		\leq \frac{|Z| \|f''\|_\infty}{2 \sqrt n} + \frac{m_0 \|f'''\|_\infty}{3 \sqrt  n}.
\end{equation}
Let us define the \emph{Stein operator} $\bb S_\nu$ as
\[
\bb S_\nu f(x) := \nu^2 f'(x) -x f(x)
\]
for every $f \in \mc C_b^1(\bb R)$ and every $x \in \bb R$. We see that for $\nu^2 = \frac{1}{2} \gamma(\frac{\ell^2}{n^2})$, \eqref{rancagua} can be rewritten as
\begin{equation}
\label{graneros}
\big| L_{n,\ell}^{\mathrm{ex}} f(Z) - \bb S_\nu f'(Z) \big| \leq \frac{|Z| \|f''\|_\infty}{2 \sqrt n} + \frac{m_0 \|f'''\|_\infty}{3 \sqrt n}.
\end{equation}
At this point, we can use \emph{Stein's lemma} in order to obtain a quantitative CLT for the random variable $Z$:

\begin{proposition}[Stein's Lemma]
\label{stein}
Let $W \sim \mc N(0,\nu^2)$ for $\nu>0$ and let $h \in \mc C_b^1(\bb R)$. Let $f_h$ be the only bounded solution of
\[
-x f (x) + \nu^2 f'(x) = h(x) -\bb E[h(W)].
\]
We have that
\[
\|f_h\|_\infty \leq 2 \|h' \|_\infty, \quad \|f_h'\|_\infty \leq \sqrt{\frac{2}{\pi \nu^2}}\|h'\|_\infty \text{ and } \|f_h''\|_\infty \leq \frac{2}{\nu^2} \|h'\|_\infty.
\]
\end{proposition}

In the case $\nu=1$, this lemma is stated without proof as Lemma 2.5 in \cite{Ros} and stated and proved as Lemma 2.4 in \cite{CheGolSha}. The case $\nu \neq 1$ follows from a simple scaling argument.

Taking $f(x) := \int_{-\infty}^x f_h(y) dy$ and using the fundamental relation $ \int  L_{n,\ell}^{\mathrm{ex}} f(Z) du(\sigma) =0$, we see that
\[
\big| \bb E_n^{u(\sigma)} [ h(Z)] - \bb E[h(W)] \big| \leq \frac{C}{\sqrt n} \Big( \frac{m_0}{\nu^2} + \frac{\bb E_n^{u(\sigma)}[|Z|]}{\nu}\Big) \|h'\|_\infty
\]
for every $h \in \mc C_b^1(\bb R)$. In order to conclude, we just need to estimate $\bb E_n^{u(\sigma)}[|Z|]$.
From Lemma \ref{l3.5.1} and Cauchy-Schwartz inequality,
\[
\bb E_n^{u(\sigma)}[|Z|] \leq \sqrt{\tfrac{n}{n-1}} m_0(1-m_0).
\]
Observe as well that
\[
\nu^2 = \tfrac{1}{2} \gamma \big( \tfrac{\ell^2}{n^2} \big) = m_0^2(1-m_0)^2.
\]
We conclude that 
\[
\big| \bb E_n^{u(\sigma)} [ h(Z)] - \bb E[h(W)] \big| \leq
		\frac{C}{\sqrt n} \Big( 1+ \frac{1}{m_0}\Big) \|h'\|_\infty,
\]
from where we conclude the following estimate:

\begin{lemma}
\label{l3.8} There exists a finite constant $C$ such that
\[
d_K(\zeta_0^1, \mc N(0,\nu^2)) \leq \frac{C}{m_0(1-m_0) \sqrt n}
\]
for every $n \in \bb N$ and every $\sigma \in \Onast$, where 
\begin{equation}
\label{rengo}
\nu = m_0(1-m_0)
\end{equation}
and $m_0:= M(\sigma)$.
\end{lemma}

\subsection{The quantitative CLT under $\bb P_n^{u(\sigma)}$}

Let $(\veczt; t \geq 0)$ be the solution of
\begin{equation}
\label{sde2}
d z_t^i = \big( a^i (\vec 1 \cdot \veczt) - \big(1+\tfrac{a+b}{n} \big) z_t^i \big) dt + \sqrt{a^i G^i \! (m_t,m_t)} d B_t^i, \quad i \in Q,
\end{equation}
with initial condition $(Z,-Z)$, where $Z$ has law $\mc N(0,\nu^2)$ and $\nu$ is defined in \eqref{rengo}. 
The aim of this section is to prove the following Quantitative CLT for $\veczetat$:

\begin{lemma}
\label{l3.9}
There exists a finite constant $C(a,b)$ such that
\[
d_K( \veczetat, \veczt) \leq C(a,b) \Big( \frac{t^{1/4} + t^{3/4}}{n^{1/6}} + \frac{1}{m_0(1-m_0) \sqrt n}\Big)
\]
for every $n \in \bb N$, every $\sigma \in \Onast$ and every $t \geq 0$.
\end{lemma}

\begin{proof}
Let $f \in \mc C_b(\bb R^Q; \bb R)$ and $t \geq 0$, and let $(g_{s,t}; s \in [0,t])$ be the solution of \eqref{fokpla3} with final condition $g_{t,t} = f$. Observe that
\[
\bb E[ f(\veczt)] = \bb E[ g_{0,t}(\veczzero)].
\]
By Dynkin's formula, Lemma \ref{l3.6}, Lemma \ref{l3.6.1} and Lemma \ref{l3.6.11}
\[
\begin{split}
\big| \bb E_n^{u(\sigma)}[f(\veczetat)] - \bb E[g_{0,t}(\veczetazero)] \big| 
		&\leq C(a,b) \int_0^t \Big( \frac{\|D^2 f\|_\infty \|\veczetas\|_{\ell^1}}{\sqrt n} + \frac{\|D^3 f\|_\infty}{\sqrt n} + \frac{\| D^4 f\|_\infty}{n} \Big)ds\\
		&\leq C(a,b) \Big( \frac{\|D^2 f\|_\infty(t+ t^{3/2})}{\sqrt n} + \frac{\|D^3f\|_\infty t}{\sqrt n} + \frac{\|D^4 f\|_\infty t}{n} \Big).
\end{split}
\]
Recall that we want to compare $\bb E_n^{u(\sigma)}[f(\veczetat)]$ with $\bb E[f(\veczt)]=  \bb E[ g_{0,t}(\veczzero)]$. Using Lemma \ref{l3.8} and Lemma \ref{l3.6.11},
\[
\begin{split}
\big| \bb E[g_{0,t}(\veczetazero)] -  \bb E[ g_{0,t}(\veczzero)] \big| 
		&\leq \|D g_{0,t}\|_\infty d_K( \veczetazero, \veczzero)\\
		&\leq \frac{C \|D f\|_\infty }{m_0(1-m_0) \sqrt n}.
\end{split}
\]
We conclude that
\[
\begin{split}
\big| \bb E_n^{u(\sigma)}[f(\veczetat)] - \bb E[f(\veczt)] \big| 
		&\leq C(a,b) \Big( \frac{ \|D f\|_\infty }{m_0(1-m_0) \sqrt n}+\frac{\|D^2 f\|_\infty(t+ t^{3/2})}{\sqrt n} + \\
		& \quad \quad \quad \quad \quad \quad + \frac{\|D^3f\|_\infty t}{\sqrt n} + \frac{\|D^4 f\|_\infty t}{n} \Big).
\end{split}
\]
Using Lemma \ref{reglem}, and keeping only the dominant terms in $n$ and $t$, the lemma is proved.
\end{proof}

\subsection{Asymptotic behavior of the associated diffusions}

Lemmas \ref{l3.5} and \ref{l3.9} reduce the computation of the distance between $\vecxit$ and $\veczetat$ to the computation of the distance between the Gaussian variables $\vecyt$ and $\veczt$. The aim of this section is to estimate this distance by means of the derivation of the asymptotic behavior of the variables $\vecyt$, $\veczt$ for $1 \ll t \ll n^\gamma$, with $\gamma \in (0,1)$ chosen properly. We will prove the following estimate:

\begin{lemma}
\label{asympt}
Let $W,W'$ be independent, standard Gaussian random variables and define $\bar{z}_t =(\bar{z}_{\!t}^i; i \in Q)$ as
\[
\bar{z}_{\!t}^i := (2i-1) \sqrt{\tfrac{1}{2} m_0(1-m_0) G(m_0)} W + a^i \sqrt{G(m_0)t} W'
\]
for each $i \in Q$. There exists C(a,b) finite such that
\[
d_K(\vecyt, \bar{z}_t) + d_K(\veczt, \bar{z}_t) \leq C(a,b) \Big( \frac{t^{3/2}}{n \sqrt{m_0(1-m_0)}} + \frac{t e^{-t}}{\sqrt{m_0(1-m_0)}} \Big)
\]
for every $t \geq 1$, every $n \in \bb N$ and every $\sigma \in \Onast$.
\end{lemma}

First we will describe the asymptotic behavior of $\vecyt$, the behavior of $\veczt$ can be derived in the same way. It will be convenient to diagonalize the drift in equation \eqref{sde}. This is accomplished defining the variables
\[
y_t := \vec{1} \cdot \vecyt, \quad \ytast := y_{\!t}^1- a^1 \! y_t.
\]
Observe that
\[
d y_t = -\big( \tfrac{a+b}{n} \big) y_t dt + \sumiq \sqrt{a^i G^i \! (\vecmt)} d B_t^i,
\]
and
\[
d \ytast = - \big(1+\tfrac{a+b}{n} \big) \ytast dt + m_0(1-m_0) \sumiq (2i-1) \sqrt{\frac{G^i \! (\vecmt)}{a^i}} d B_t^i.
\]
By the variation of parameters formula, 
\begin{equation}
\label{teno}
y_t = \int_0^t e^{-(\frac{a+b}{n})(t-s)} \sumiq \sqrt{a^i G^i \! (\vecms)} d B_s^i,
\end{equation}
and
\begin{equation}
\label{losquenes}
\ytast = \int_0^t e^{-(1+\frac{a+b}{n})(t-s)} m_0(1-m_0) \sumiq (2i-1)\sqrt{\frac{G^i \! (\vecms)}{a^i}} d B_s^i.
\end{equation}
Before we enter into the derivation of the asymptotic behavior of $y_t$ and $\ytast$, we need some estimates on the asymptotic behavior of $\vecmt$:

\begin{lemma}
\label{l3.14} There exists a constant $C(a,b)$ such that:
\begin{itemize}
\item[a)] for every solution $(m_t; t \geq 0)$ of \eqref{EDO} and every $t \geq 0$,
\begin{equation}
\label{coquimbo2}
|m_t - m_0| \leq \frac{C(a,b) t}{n},
\end{equation}
\item[b)]  for every solution $(m_t; t \geq 0)$ of \eqref{EDO} and every $t \geq 0$,
\begin{equation}
\label{pumalquin}
\frac{1}{m_t(1-m_t)} \leq \frac{C(a,b)}{m_0(1-m_0)},
\end{equation}
\item[c)] for every solution $(\vecmt; t \geq 0)$ of \eqref{EDO2}, every $t \geq 0$ and every $i \in Q$,
\begin{equation}
\label{romeral}
|m_t^i-m_t| \leq (1-a^i) e^{-t}.
\end{equation}
\end{itemize}
\end{lemma}

\begin{proof}
Estimate \eqref{coquimbo2} is just a restatement of \eqref{coquimbo}. In order to prove \eqref{pumalquin}, observe that $m_t$ converges, as $t \to \infty$, to $\frac{a}{a+b}$. Therefore, if $m_0(1-m_0) \leq \frac{ab}{(a+b)^2}$, then $m_t(1-m_t)$ is minimal at $t =0$. In order to prove \eqref{romeral}, it is enough to observe that
\[
\tfrac{d}{dt} (m_t^i-m_t) = -\big(1+\tfrac{a+b}{n} \big) (m_t^i-m_t)
\]
and that $|m_0^i-m_0| = 1-a^i$. 
\end{proof}

Let us describe the asymptotic behavior of $y_t$ when $t \to \infty$. Observe that
\[
\sumiq a^i G^i \! (\vecmt) = G(m_t).
\]
Let $(\bar{B}_t; t \geq 0)$ be given by
\begin{equation}
\label{lolol}
\bar{B}_t :=  \sumiq \sqrt{a^i} B_t^i
\end{equation}
for every $t \geq 0$. Observe that $(\bar{B}_t; t \geq 0)$ is a standard Brownian motion. Let $\omega_\infty$ be given by
\[
\omega_\infty := \int_0^\infty \sumiq \big( \sqrt{ a^i G^i \! (\smash{\vecms})} - \sqrt{ a^i G(\smash{m_s})} \big) d B_s^i.
\]
We have the following estimate:

\begin{lemma}
\label{l3.15} There exists $C(a,b)$ finite such that for every $t \geq 0$, every $n \in \bb N$ and every $\sigma \in \Onast$,
\[
\Var \big( y_t -\omega_\infty- \sqrt{G(m_0)} \bar{B}_t \big) \leq C(a,b)\Big(\frac{t^3}{n^2 m_0(1-m_0)} + e^{-2t}\Big).
\]
\end{lemma}

\begin{proof}
The idea is to replace \eqref{teno} by simpler expressions, until we end up with the variables $\omega_\infty$ and $\bar{B}_t$. First we remove the exponential factor in \eqref{teno}. We have that
\[
\begin{split}
\Var \Big( y_t - \int_0^t \sumiq \sqrt{a^i \smash{G^i \! (\vecms)}} d B_s^i \Big) 
		&= \int_0^t \Big(1- e^{-(\frac{a+b}{n})(t-s)}\Big)^2 \sumiq a^i G^i \! (\vecms) ds\\
		&\leq \| \vecG\|_\infty \int_0^t \frac{(a+b)^2 (t-s)^2}{n^2} ds \leq \frac{C(a,b) t^3}{n^2}.
\end{split}
\]
Therefore, it is enough to consider 
\[
\int_0^t \sumiq  \sqrt{a^i \smash{G^i \! (\vecms)}} d B_s^i.
\]
Observe that
\[
\int_0^t \sumiq  \sqrt{a^i \smash{G^i \! (\vecms)}} d B_s^i - \int_0^t \sumiq  \sqrt{a^i \smash{G(m_s)}} d B_s^i 
		= \omega_\infty - \int_t^\infty \sumiq  \big( \sqrt{a^i \smash{G^i \! (\vecms)}} - \sqrt{a^i G(\smash{m_s})} \big) d B_s^i.
\]
Using the elementary inequality 
\begin{equation}
\label{molina}
(\sqrt{u}-\sqrt{v})^2 \leq \frac{(u-v)^2}{u},
\end{equation}
we see that
\[
\begin{split}
\Var \Big( \int_t^\infty \sumiq  \big( \sqrt{a^i \smash{G^i \! (\vecms)}} 
		&- \sqrt{a^i G(\smash{m_s})} \big) d B_s^i \Big)
		= \int_t^\infty \sumiq  \big( \sqrt{a^i \smash{G^i \! (\vecms)}} - \sqrt{a^i G(\smash{m_s})} \big)^2 ds\\
		&\leq \int_t^\infty \sumiq \frac{a^i (G^i \! (\vecms)-G(m_s))^2}{G(m_s)}ds
		\leq \int_t^\infty \sumiq \frac{C(a,b) a^i \! (1-a^i)^2 e^{-2s}}{m_0(1-m_0)} \\
		&\leq C(a,b) e^{-2t}.
\end{split}
\]
Therefore, we are left to consider the integral
\[
\int_0^t \sumiq \sqrt{a^i G(\smash{m_s})} d B_s^i.
\]
In order to finish the proof of the lemma, we need to replace $G(m_s)$ with $G(m_0)$. We have that
\[
\begin{split}
\Var \Big( \int_0^t \sumiq \big(\sqrt{ a^i G(\smash{m_s})} 
		&-\sqrt{a^ i G(\smash{m_0})} \big) d B_s^i \Big) 
		= \int_0^t  \big(\sqrt{ G(\smash{m_s})} -\sqrt{ G(\smash{m_0})} \big)^2 ds\\
		&\quad \quad \leq \int_0^t  \frac{(G(m_s)-G(m_0))^2}{G(m_0)} ds \leq \int_0^t \frac{C(a,b) (m_s-m_0)^2}{m_0(1-m_0)} ds\\
		&\quad \quad\leq \frac{C(a,b) t^3}{n^2 m_0(1-m_0)}.
\end{split}
\]
Since
\[
\int_0^t \sqrt{ G(m_0)} d \bar{B}_s = \sqrt{G(m_0)} \bar{B}_t,
\]
the lemma is proved.
\end{proof}

The random variables $\omega_\infty$ and $\bar{B}_t$ are correlated in a non-trivial way, so Lemma \ref{l3.15} does not give a complete picture of the asymptotic behavior of $y_t$. However, since $y_t$ is a Gaussian random variable, its asymptotic behavior is completely characterized by its variance. We have the following estimate:

\begin{lemma}
\label{l3.15.1}
There exists $C(a,b)$ finite such that
\[
\big| \Var(y_t) - G(m_0) t \big| \leq \frac{C(a,b) t^2}{n}
\]
for every $n \in \bb N$, every $\sigma \in \Onast$ and every $t \geq 0$.
\end{lemma}

\begin{proof}
Observe that
\[
\begin{split}
\Var(y_t) 
		&= \int_0^t e^{-2 (\frac{a+b}{n})(t-s)} \sumiq a^i G^i (\vecms) ds = \int_0^t e^{-2 (\frac{a+b}{n})(t-s)} G (m_s) ds
\end{split}
\]
and therefore
\[
\big| \Var(y_t) - G(m_0) t \big| 
		\leq \int_0^t \Big| e^{-2 (\frac{a+b}{n})(t-s)} -1 \Big| G(m_s)ds + \int_0^t \big| G(m_s) - G(m_0)\big| ds.
\]
The lemma follows from the estimates of Lemma \ref{l3.14}.
\end{proof}

\begin{remark}
Since $y_t$ is a Gaussian process, the asymptotic behavior of $y_t$ for $1 \ll t \ll \sqrt n$ can be described by Lemma \ref{l3.15.1}, without appealing to Lemma \ref{l3.15}. The point of Lemma \ref{l3.15} is to be able to describe the asymptotics of the \emph{joint law} of $(y_t,\ytast)$, for which we need to establish the correlations between these two processes. We anticipate that $y_t$ and $\ytast$ are asymptotically independent.
\end{remark}

In order to describe the asymptotic behavior of $\ytast$, let us extend the definition of  $(B_t^i; t \geq 0, i \in Q)$ to a two-sided Brownian motion $(B_t^i; t \in \bb R, i \in Q)$.
Define
\[
\omega_{t,\infty}^\ast := \int_{-\infty}^t e^{-(t-s)} m_0(1-m_0) \sumiq (2i-1) \sqrt{\frac{G(m_0)}{a^i}} d B_{s}^i.
\]
Observe that 
\[
\Cov(\omega_{t,\infty}^\ast, \bar{B}_t)= \int_0^t e^{-(t-s)} m_0(1-m_0) \sumiq (2i-1) \sqrt{\frac{G(m_0)}{a^i}} \sqrt{a^i} ds =0.
\]
We have the following result:

\begin{lemma}
\label{l3.16}
There exists $C(a,b)$ finite such that for every $t \geq 0$, every $n \in \bb N$ and every $\sigma \in \Onast$,
\[
\Var( \ytast -\omega_{t,\infty}^\ast) \leq C(a,b) \Big( \frac{t^2}{n^2} + (1+ t) e^{-2t} \Big).
\]
\end{lemma}

\begin{proof}
As above, the idea is to repeatedly replace $\ytast$ by simpler expressions and to estimate the variance of the different error terms until we end up with $\omega_\infty^\ast$. Recall the representation \eqref{losquenes} for $\ytast$. We have that
\[
\begin{split}
\Var \Big( \ytast - \int_0^t e^{-(t-s)} m_0(1-m_0) 
		&\sumiq (2i-1) \sqrt{\frac{G^i \! (\vecms)}{a^i}} d{B}_{s}^i \Big) \\
		&\!\!= \int_0^t e^{-2(t-s)} \Big( 1- e^{-(\frac{a+b}{n})(t-s)}\Big)^2 m_0^2(1-m_0)^2 \sumiq \frac{G^i \! (\vecms)}{a^i} ds \\
		&\!\!\leq \frac{C(a,b)m_0^2(1-m_0)^2}{n^2} \int_0^t e^{-2(t-s)} s^2 ds \leq \frac{C(a,b) m_0(1-m_0)t^2}{n^2}.
\end{split}
\]
Therefore, we are left to consider
\[
\int_0^t e^{-(t-s)} m_0(1-m_0) \sumiq (2i-1) \sqrt{\frac{G^i \! (\vecms)}{a^i}} d{B}_{s}^i.
\]
Now we will replace $G^i \! (\vecms)$ by $G(m_{s})$. Using \eqref{molina}, \eqref{pumalquin} and \eqref{romeral}, we have that
\[
\begin{split}
\Var \Big( \int_0^t e^{-(t-s)} m_0(1-m_0) 
		&\sumiq (2i-1) \Big( \sqrt{\frac{G^i \! (\vecms)}{a^i}} - \sqrt{\frac{G(m_{s})}{a^i}}\Big) d{B}_{s}^i \Big)\\
		&= \int_0^t e^{-2(t-s)} m_0^2(1-m_0)^2 \sumiq \frac{1}{a^i} \big( \sqrt{G^i \! (\smash{\vecms})} - \sqrt{G(m_{s})}\big)^2 ds\\
		&\leq C(a,b) \int_0^t e^{-2(t-s)} \sumiq \frac{m_0^2(1-m_0)^2(1-a^i)^2 e^{-2s}}{a^i m_0(1-m_0)} ds\\
		&\leq C(a,b) te^{-2t}.
\end{split}
\]
Therefore, it is enough to consider 
\[
\int_0^t e^{-(t-s)} m_0(1-m_0) \sumiq (2i-1) \sqrt{\frac{G(m_{s})}{a^i}} d {B}_{s}^i.
\]
Now we will replace $m_{s}$ by $m_0$:
\[
\begin{split}
\Var \Big( \int_0^t e^{-(t-s)} m_0(1-m_0) 
		&\sumiq (2i-1) \Big( \sqrt{\frac{G(m_{s})}{a^i}} - \sqrt{\frac{G(m_0)}{a^i}} \Big)d{B}_{s}^i \Big)\\
		&= \int_0^t e^{-2(t-s)} m_0^2(1-m_0)^2 \sumiq \frac{1}{a^i} \big( \sqrt{G(m_{s})}-\sqrt{G(m_0)}\big)^2 ds\\
		&\leq C(a,b) \int_0^t e^{-2(t-s)} m_0^2 (1-m_0)^2 \sumiq \frac{s^2}{a^i m_0(1-m_0) n^2} ds\\
		&\leq \frac{C(a,b) t^2}{n^2}.
\end{split}
\]
Finally, we observe that
\[
\begin{split}
\Var \Big( \int_t^\infty e^{-s} m_0(1-m_0) 
		&\sumiq (2i-1) \sqrt{\frac{G(m_0)}{a^i}} d {B}_s^i \Big) \\
		&= \int_t^\infty e^{-2s} m_0^2(1-m_0)^2 \sumiq \frac{G(m_0)}{a^i} ds
		\leq C(a,b) m_0(1-m_0) e^{-2t}.
\end{split}
\]
Collecting all the estimates, the lemma is proved.
\end{proof}

In order to fully describe the asymptotic behavior of $(y_t, \ytast)$, we are left to show that $\omega_\infty$ and $\omega_{t,\infty}^\ast$ are asymptotically independent. It will be more convenient to directly estimate the covariance between $y_t$ and $\omega_{t,\infty}^\ast$. We have the following estimate:

\begin{lemma}
\label{l3.16.1}
There exists $C(a,b)$ finite such that for every $t \geq 0$, every $n \in \bb N$ and every $\sigma \in \Onast$,
\[
\big|\Cov( y_t, \omega_{t,\infty}^\ast)\big| \leq C(a,b) \sqrt{m_0(1-m_0)} te^{-t}.
\]
\end{lemma}

\begin{proof}
We have that
\[
\begin{split}
\Cov(y_t, \omega_{t,\infty}^\ast) 
		&= \int_0^t e^{-(1+\frac{a+b}{n}) (t-s)} m_0(1-m_0) \sumiq (2i-1) \sqrt{G^i \! (\smash{\vecms})G(m_0)} ds \\
		&= \int_0^t e^{-(1+\frac{a+b}{n}) (t-s)} m_0(1-m_0) \sqrt{G(m_0)} \big( \sqrt{G^1(\vecms)} - \sqrt{G^0(\vecms)}\big) ds
\end{split}
\]
Recall \eqref{romeral}. Using \eqref{molina}, we see that
\[
\begin{split}
\big| \Cov(y_t, \omega_{t,\infty}^\ast) \big| 
		&\leq C(a,b) m_0(1-m_0) \sqrt{G(m_0)}\int_0^t e^{-t+s} \sumiq \frac{(1-a^i) e^{-s}}{\sqrt{G^1(\vecms)}} ds\\
		&\leq C(a,b) \sqrt{m_0(1-m_0)} t e^{-t},
\end{split}
\]
which proves the lemma.
\end{proof}

It turns out that the asymptotic behavior of $\veczt$ is identical to the behavior of $\vecyt$. Defining the variables $z_t := \vec{1} \cdot \veczt$, $\ztast := z_t^1-a^1 z_t$, we see that
\[
z_t = \int_0^t e^{-(\frac{a+b}{n})(t-s)} \sumiq \sqrt{a^i G(m_s)} d B_s^i
\]
and
\[
\ztast = Z e^{-(1+\frac{a+b}{n})t} + \int_0^t e^{-(1+\frac{a+b}{n})(t-s)} m_0 (1-m_0) \sumiq (2i-1) \sqrt{\frac{G(m_s)}{a^i}} d B_s^i.
\]
Observe that, differently from $y_t$ and $\ytast$, $z_t$ and $\ztast$ are independent. Recall definition \eqref{lolol}. 
We have the following estimates:

\begin{lemma}
\label{l3.17} For every $t \geq 0$, every $n \in \bb N$ and every $\sigma \in \Onast$, 
\[
\Var \big( z_t - \sqrt{G(m_0)} \bar{B}_t \big) \leq \frac{C(a,b) t^3}{n^2 m_0(1-m_0)}
\]
and
\[
\Var \big( \ztast - \omega_{t,\infty}^\ast) \leq C(a,b) \Big( m_0(1-m_0) e^{-2t} +\frac{t^2}{n^2} \Big)
\]
\end{lemma}

The proof of these estimates is almost exactly equal to the proofs of Lemma \ref{l3.15} and \ref{l3.16}, and therefore we omit it. The only difference is the appearance of the random variable $Z$. The bound $\Var(Z) \leq 2m_0^2(1-m_0)^2$ leads to the term $m_0(1-m_0)$ in the second estimate.

Now we observe that $y_t^1 = a^1y_t + \ytast$, $y_t^0 = a^0 y_t - \ytast$. In order to describe the asymptotic behavior of $\vecyt$ and $\veczt$, we can forget about the exact representation of the random variables $\bar{B}_t$, $\omega_\infty$ and $\omega_{t,\infty}^\ast$, and simply recall that $(\bar{B}_t, \omega_\infty, \omega_{t,\infty}^\ast)$ is a centered Gaussian vector. Observe that
\[
\Var(\omega_{t,\infty}^\ast) = \int_0^\infty e^{-2s} m_0^2(1-m_0)^2 \sumiq \frac{G(m_0)}{a^i} ds = \tfrac{1}{2} m_0(1-m_0) G(m_0). 
\]

In order to compare $\vecyt$ and $\veczt$, we will use the following lemma:

\begin{lemma}
\label{gaussiancoupling}
Let $(X,Y,Z)$ be a Gaussian vector. Assume that $X$, $Z$ are independent, and that $Z \neq 0$. There exists a random variable $\tilde{Y}$ such that $(Y,Z)$ and $(\tilde{Y},Z)$ have the same law and such that
\[
\bb E[(X - \tilde Y)^2] \leq \frac{2(\Var(X)-\Var(Y))^2}{\Var(X)} + \frac{3 \Cov (Y,Z)^2}{\Var(Z)}.
\]
\end{lemma}

\begin{proof}
Let us assume that $\tilde Y = \alpha X + \beta Z$, for some $\alpha, \beta \in \bb R$. Without loss of generality, we can assume $\alpha >0$. We have that 
\[
\Var(\tilde Y) = \alpha^2 \Var(X) + \beta^2 \Var(Z), \quad \Cov(\tilde Y,Z) = \beta \Var(Z). 
\]
Therefore, $(\tilde Y,Z)$ and $(Y,Z)$ have the same law if $\alpha$, $\beta$ satisfy
\[
\beta = \frac{\Cov(Y,Z)}{\Var(Z)} \text{ and } \alpha^2 = \frac{\Var(Y)}{\Var(X)} - \frac{\Cov(Y,Z)^2}{\Var(X) \Var(Z)}.
\]
Observe that the Cauchy-Schwartz estimate
\begin{equation}
\label{talca}
\big| \Cov(Y,Z) \big| \leq \Var(Y)^{1/2} \Var(Z)^{1/2}
\end{equation}
implies that $\alpha$ is well defined. For these values of $\alpha$, $\beta$, 
\[
\begin{split}
\bb E[(X-\tilde Y)^2] 
		&= (\alpha-1)^2\Var(X) + \beta^2 \Var(Z) \\
		&= \Big( \sqrt{ \Var(Y) - \frac{\Cov(Y,Z)}{\Var(Z)}}-\sqrt{\Var(X)}\Big)^2 + \frac{\Cov(Y,Z)^2}{\Var(Z)}
\end{split}
\]
Using \eqref{molina} and \eqref{talca}, we see that
\[
\begin{split}
\bb E[(X-\tilde Y)^2] 
		&\leq \frac{1}{\Var(X)} \Big( \Var(Y) - \Var(X) - \frac{\Cov(Y,Z)^2}{\Var(Z)} \Big)^2 +  \frac{\Cov(Y,Z)^2}{\Var(Z)}\\
		&\leq \frac{2 (\Var(Y)-\Var(X) \big)^2}{\Var(X)} + \frac{2 \Cov(Y,Z)^4}{\Var(X) \Var(Z)^2} + \frac{\Cov(Y,Z)^2}{\Var(Z)}\\
		&\leq \frac{2 (\Var(Y)-\Var(X) \big)^2}{\Var(X)} + \frac{3 \Cov(Y,Z)^2}{ \Var(Z)},
\end{split}
\]
as we wanted to show.
\end{proof}

\begin{proof}[Proof of Lemma \ref{asympt}] Recall that in the hypothesis of Lemma \ref{asympt}, $t \geq 1$. Observe that $\bar{z}$ has the same law of
\[
\big((2i-1)\omega_{t,\infty}^\ast+ a^i \sqrt{G(m_0)} \bar{B}_t; i \in Q \big).
\]
Let us define $\bar{y}_t$ as
\[
\bar{y}_{\!t}^i := (2i-1) \omega_{t,\infty}^\ast + a^i y_t
\]
for every $i \in Q$. By \eqref{dual}, If $X$, $Y$ are square-integrable random variables, then
\[
d_K(X,Y)^2 \leq \bb E\big[ \|X-Y\|^2 \big].
\]
Therefore, by Lemma \ref{gaussiancoupling}, with $X = \sqrt{G(m_0)} \bar{B}_t$, $Y = y_t$ and $Z = \omega_{t,\infty}^\ast$,
\[
d_K(\bar{y}_t, \bar{z}_t)^2 
		\leq \sumiq (a^i)^2 \bb E[(X - \tilde Y)]
		\leq \sumiq (a^i)^2 \Bigg( \frac{2( G(m_0) t - \Var(y_t))^2 }{G(m_0)t} + \frac{3\Cov(y_t,\omega_{t,\infty}^\ast)^2}{\Var(\omega_{t,\infty}^\ast)}\Bigg).
\]
Observe that $(a^1)^2+(a^0)^2 \leq \frac{1}{2}$. Using Lemmas \ref{l3.15.1} and \ref{l3.16.1} we conclude that
\[
\begin{split}
d_K(\bar{y}_t, \bar{z}_t)^2 
		&\leq C(a,b) \Big( \frac{t^3}{n^2 m_0(1-m_0)} + \frac{m_0(1-m_0)t^2e^{-2t}}{m_0(1-m_0)G(m_0)}\Big)\\
		&\leq \frac{C(a,b)}{m_0(1-m_0)}\Big( \frac{t^3}{n^2} + t^2e^{-2t}\Big).
\end{split}
\]
This estimate is compatible with the lemma. Recall that
\[
y_{\!t}^i = (2i-1) y_t^\ast +a^i y_t.
\]
Therefore,
\[
d_K(\vecyt, \bar{y}_t)^2 \leq 2 \Var(y_t^\ast - \omega_{t,\infty}) \leq C(a,b)\Big( \frac{t^2}{n^2} + (1+t)e^{-2t}\Big) 
\]
This proves the lemma for $\vecyt$. In order to prove the lemma for $\veczt$, we need to compare the vectors
\[
\big( (2i-1) \omega_{t,\infty}^\ast + a^i \sqrt{G(m_0)} \bar{B}_t; i \in Q \big),
\]
\[
\big( (2i-1) z_t^\ast + a^i z_t; i \in Q \big).
\]
Using Lemma \ref{l3.17}, we see that
\[
\begin{split}
d_K(\veczt, \bar{z}_t)^2 
		&\leq  4\Var \big(z_t -\sqrt{G(m_0)} \bar{B}_t \big) + \Var \big( z_t^\ast - \omega_{t,\infty}^\ast \big)\\
		&\leq C(a,b) \Big( \frac{t^3}{n^2 m_0(1-m_0)} + m_0(1-m_0) e^{-2t} \Big),
\end{split}
\]
which proves the lemma for $\veczt$.
\end{proof}

\subsection{Proof of Theorem \ref{phase1}}

In this section we prove Theorem \ref{phase1}. According to Lemma \ref{acopl}, it is enough to consider the local densities. Let us denote by $\vecYt$ the law of $\vecX(\sigma)$ with respect to $\pens$ and by $\vecZt$ the law of $\vecX(\sigma)$ with respect to $\penus$. We have that
\[
Y_t^i = a^i \! n m_t^i + \sqrt{n} \xi_t^i, \quad Z_t^i  = a^i \! n m_t + \sqrt{n} \zeta_t^i.
\]
Recall the definition of $\bar{z}_t$ given in Lemma \ref{asympt} and define
\[
\bar{Y}_t^i := a^i \! n m_t^i + \sqrt{n} \bar{z}_t^i, \quad \bar{Z}_t^i := a^i \! n m_t + \sqrt{n} \bar{z}_t^i.
\]
According to Lemmas \ref{l3.5}, \ref{l3.9} and \ref{asympt}, 
\[
\begin{split}
\frac{1}{\sqrt n} \big| d_K(\vecYt, \vecZt) - d_K(\bar{Y}_t, \bar{Z}_t) \big| 
		&\leq d_K(\vecxit, \vecyt) + d_K(\vecyt,\bar{z}_t) + d_K(\bar{z}_t, \veczt) + d_K(\veczt, \veczetat)\\
		&\leq C(a,b) \Bigg( \frac{t^{1/4} + t^{3/4}}{n^{1/6}} + \frac{t^{3/2}}{n \sqrt{m_0(1-m_0)}} \\
		&\quad \quad \quad \quad+ \frac{te^{-t}}{\sqrt{m_0(1-m_0)}} + \frac{1}{m_0(1-m_0) \sqrt{n}}\Bigg).
\end{split}
\]
We will see that for our choices of times $t$ and initial conditions $m_0$, the right-hand side of this estimate vanishes as $n \to \infty$. 
Therefore, it will be enough to analyze the distance $d_K(\bar{Y}_t, \bar{Z}_t)$. But $\bar{Y}_t$ is a translation of $\bar{Z}_t$. From Lemma \ref{la3}, we see that
\[
\frac{1}{\sqrt n} d_K(\bar{Y}_t, \bar{Z}_t) = \sqrt{n} \sumiq a^i |m_t^i-m_t| = 2 \sqrt{n} m_0(1-m_0) e^{-(1+\frac{a+b}{n})t}.
\]
In particular, for every $\tau \in \bb R$, taking $t_n := \frac{1}{2} \log n + \log m_0(1-m_0) + \tau$, we have that
\[
\lim_{n \to \infty} \frac{1}{\sqrt n} d_K(\bar{Y}_{t_n}, \bar{Z}_{t_n} ) = 2 e^{-\tau}.
\]
Recall that $m_0 \in \big[\frac{1}{n}, 1-\frac{1}{n} \big]$ and therefore $t_n \leq \frac{3}{2} \log n + \tau$. For this choice of times $t_n$, we see that
\[
\frac{1}{\sqrt n} \big| d_K(\vecYtn, \vecZtn) - d_K(\bar{Y}_{t_n}, \bar{Z}_{t_n}) \big| 
		\leq C(a,b) \Big( \frac{(\log n)^{3/4}}{n^{1/6}} + \frac{\log n }{ n^{1/2- 3\alpha/2}}\Big),
\]
Since $\alpha < \frac{1}{3}$, 
\[
\lim_{n \to \infty} \frac{1}{\sqrt n} \big| d_K(\vecYtn, \vecZtn) - d_K(\bar{Y}_{t_n}, \bar{Z}_{t_n}) \big| =0
\]
for every $\tau \in \bb R$, which proves the theorem.

\appendix

\section{}

\begin{lemma}[Regularization Lemma] 
\label{reglem}
Let $\mu$, $\nu$ be two probability measures in $\mc P_1(\bb R^d)$ and let $\ell \geq 2$. Assume that the exist constants $(a_i; i =1,\dots,\ell)$ such that
\[
\Big| \int \! f d \mu - \int \! f d \nu \Big| \leq a_i \|D^i f\|_\infty
\]
for every $f \in \mc C_b^\infty(\bb R^d; \bb R)$. There exists a constant $C= C(d,\ell)$ such that
\[
d_K(\mu,\nu) \leq C \max_i \{a_i^{1/i} \}.
\]
\end{lemma}

\begin{proof}
Let $g$ be the density of a standard Gaussian in $\bb R^d$:
\[
g(x) := \frac{1}{(2\pi)^{d/2}} e^{-\|x\|^2/2} \; \forall x \in \bb R^d.
\]
For $\theta >0$, let $g_\theta$ be given by $g_\theta(x) := \theta^{-d} g(\theta^{-1}x)$ for all $x \in \bb R^d$. Take $f \in \mc C_b^\infty(\bb R^d; \bb R)$ and let $f_\theta := f \ast g_\theta$. Observe that
\[
\Big| \int (f-f_\theta) d \mu \Big| \leq \| D f \|_\infty \int \|x\| g_\theta(x) dx \leq C(d) \theta \|Df\|_\infty.
\]
The same estimate holds replacing $\mu$ by $\nu$. Observe as well that for every $i \geq 1$, 
\[
\|D^i f_\theta\|_\infty \leq \|D f\|_\infty \| D^{i-1} g_\theta\|_{L^1} \leq C(d,i-1) \theta^{-(i-1)} \| Df\|_\infty.
\]
We conclude that 
\[
\Big| \int f d \mu - \int f d \nu \Big| \leq \| D f\|_\infty \Big( 2C(d) \theta + \sum_{i=1}^\ell C(d,i-1) a_i \theta^{-(i-1)} \Big).
\]
Taking $\theta = \max_i a_i^{1/i}$, the lemma follows.
\end{proof}

\begin{lemma}
\label{la2}
Let $(\mc E_i, d_i)$, $i=1,2$ be complete, separable metric spaces. Let $\pi: \mc E_1 \to \mc E_2$ be a Lipschitz function such that $\llbracket f \rrbracket \leq 1$. Let $\mu, \nu \in \mc P_1(\mc E_1)$ and let $\pi_\ast \mu$, $\pi_\ast \nu$ be the pushforward of $\mu$, $\nu$ onto $\mc E_2$ respectively. We have that
\[
d_K(\pi_\ast \mu, \pi_\ast \nu) \leq d_K(\mu,\nu).
\] 
\end{lemma}
\begin{proof}
It is enough to observe that for every Lipschitz function $f: \mc E_2 \to \bb R$,
\[
\llbracket f \circ \pi \rrbracket \leq \llbracket f \rrbracket \llbracket \pi \rrbracket.
\]
\end{proof}

\begin{lemma}
\label{la3}
Let $\vecX$ be an $\bb R^d$-valued random variable, and let $\vec v \in \bb R^d$. We have that
\[
d_K(\vecX+ \vec v, \vecX) = \|\vec v\|.
\]
\end{lemma}

\begin{proof}
In one hand, the trivial coupling shows that
\[
d_K(\vecX+ \vec v, \vecX) \leq \bb E[ \|\vecX +\vec v -\vecX\| ] = \|\vec v\|.
\]
In the other hand, for $f \in \mc C_b^1(\bb R^d; \bb R)$,
\[
\big| \bb E[f(\vecX + \vec v)] -\bb E[f(\vecX)]\big| \leq \bb E [\llbracket f \rrbracket \| \vec X+\vec v - \vec X\|] = \llbracket f \rrbracket \|\vec v\|,
\]
which proves the lemma.
\end{proof}

%

\end{document}